\documentclass[a4paper,11pt,reqno]{amsart}
\usepackage{amsmath,amssymb,amsthm,bbm}
\usepackage{graphics}
\usepackage[colorlinks=true, citecolor=red, linkcolor=blue, urlcolor=black]{hyperref}
\usepackage{epsfig}
\usepackage{color}
\usepackage{enumerate}

\newcommand{\stilde}{\tilde{\sigma}}
\newcommand{\bStilde}{\tilde{\mathbb{S}}}

\numberwithin{equation}{section}

\newtheorem{theorem}{Theorem}[section]
\newtheorem{lemma}[theorem]{Lemma}
\newtheorem{proposition}[theorem]{Proposition}

\theoremstyle{definition}
\newtheorem{definition}[theorem]{Definition}
\newtheorem{remark}[theorem]{Remark}

\addtocontents{toc}{\setcounter{tocdepth}{3}}

\newcommand{\restr}{\mathop{\raisebox{-.127ex}{\reflectbox{\rotatebox[origin=br]{-90}{$\lnot$}}}}}

\newcommand{\R}{\mathbb{R}}
\newcommand{\N}{\mathbb{N}}

\newcommand{\Z}{\mathbb{Z}}
\newcommand{\bS}{\mathbb{S}}
\newcommand{\eps}{\varepsilon}
\newcommand\lt{\left}
\newcommand\rt{\right}
\def\les{\lesssim}
\def\ges{\gtrsim}

\def\diam{\operatorname{diam}}
\newcommand{\Lip}{\operatorname{Lip}}
\newcommand{\Wb}{Wb}
\newcommand{\one}{\mathbbm{1}}
\def\EE{\mathbb{E}}
\def\PP{\mathbb{P}}
\def\HH{\mathcal{H}}

\def\betinf{\underline{\mathbf{c}}^{p,d}}
\def\betsup{\overline{\mathbf{c}}_{p,d}}
\def\Sdist{\mathsf{d}_{\mathbb{S}^{d-1}}}
\def\Dzj{\Delta^z_j}
\def\tnpd{\tau_n}
\def\etnpd{\eta_n}
\def\what{\widehat}

\def\ellinf{\underline{\ell}}
\def\ellsup{\overline{\ell}}

\title[Subadditivity and optimal matching of unbounded samples] {Subadditivity and optimal matching \\  of unbounded samples
}
\author[E. Caglioti]{Emanuele Caglioti}
\address{E.C.: Dipartimento di Matematica, Guido Castelnuovo, Sapienza Universit\`a di Roma, 00185 Roma, Italy}
\email{emanuele.caglioti@uniroma1.it}
\author[M. Goldman]{Michael Goldman}
\address{M.G.:   CMAP, CNRS, \'Ecole polytechnique, Institut Polytechnique de Paris, 91120 Palaiseau,
France}
\email{michael.goldman@cnrs.fr}
\author[F. Pieroni]{Francesca Pieroni}
\address{F.P.: Dipartimento di Matematica, Guido Castelnuovo, Sapienza Universit\`a di Roma, 00185 Roma, Italy}
\email{francesca.pieroni@uniroma1.it}
\author[D. Trevisan]{Dario Trevisan}
\address{D.T.: Dipartimento di Matematica, Universit\`a degli Studi di Pisa, 56125 Pisa, Italy  }
\email{dario.trevisan@unipi.it}
\date{}
 \subjclass[2020]{60D05, 90C05, 39B62, 60F25, 49Q22}
\keywords{bipartite matching problem, optimal transport, geometric probability}

\begin{document}
\begin{abstract}
 We obtain new bounds for the optimal matching cost for empirical measures with unbounded support. For a large class of radially symmetric and rapidly decaying probability laws, we prove for the first time the asymptotic rate of convergence for the whole range of power exponents $p$ and dimensions $d$. Moreover we identify the exact prefactor when $p\le d$. We cover in particular the Gaussian case, going far beyond the currently known bounds. Our proof technique is based on approximate sub- and super-additivity bounds along a geometric decomposition adapted to some features the density, such as its radial symmetry and its decay at infinity.
\end{abstract}

\maketitle
\section{Introduction}
In this paper we investigate the behavior of the transportation cost between exponentially decaying    distributions and their empirical measure obtained from a sample of $n$ independent observations. To a large extent, we complete the program initiated in \cite{Le17,talagrand2018scaling,ledoux2019optimal,caglioti2023random}.

In order to state our main result let us introduce some notation. Letting $d\ge 1$ and $(X_i)_{i=1}^\infty$ be i.i.d.\ random variables uniformly distributed on the unit cube $Q:=(0,1)^d$ and then $\mu_n:=\sum_{i=1}^n \delta_{X_i}$ be the (non-normalized) empirical measure, the optimal matching problem concerns the study of the asymptotic behavior of
\[
 \EE\lt[W^p_p(\mu_n,n\one_{Q})\rt],
\]
where $W_p^p$ denotes the $p-$Wasserstein distance. Letting for $p,d\ge 1$,
\begin{equation}\label{defetanpd}
 \eta^{p,d}_n:=\begin{cases}
                n^{1-\frac{p}{2}} &\textrm{if } d=1,\\
                n^{1-\frac{p}{2}} (\log n)^{\frac{p}{2}} &\textrm{if } d=2,\\
                n^{1-\frac{p}{d}} &\textrm{if } d\ge 3,
               \end{cases}
\end{equation}
we define
\begin{equation}\label{defbetasup}
 \betsup:=\limsup_{n\to \infty} \frac{1}{\eta^{p,d}_n}\EE\lt[W^p_p(\mu_n,n\one_{Q})\rt].
\end{equation}
To simplify notation we write $\etnpd$ for $\eta_n^{p,d}$ when $p$ and $d$ are clear from the context.
For the case $d=1$, we refer to \cite{BoLe16} for a thorough discussion and will focus on the case $d\ge 2$ from now on. Since the work \cite{AKT84}, see also \cite[Section 4]{Le17}, it is known that for every $p\ge 1$ and every $d\ge2$, $\betsup\in(0,\infty)$.
\begin{remark}
 Let us point out that if $d\ge 3$ or $p=d=2$,  this limsup is actually a limit, see \cite{CaLuPaSi14,AmStTr16,BaBo,DeScSc13,goldman2021convergence} but this does not affect our analysis.   It is still an open problem to prove the convergence when $d=2$ and $p\neq 2$. See however \cite{AmTreVi} for recent progress in this direction.
\end{remark}
Following classical ideas in the matching problem, see e.g. \cite{yukich2006probability,BaBo,DeScSc13,ambrosio2022quadratic} we also consider a variant of the problem where mass can be transported to and from the boundary. We denote by $\Wb_Q^p$ the boundary $p-$Wasserstein distance (see \eqref{bwass} for the definition) and set
\begin{equation}\label{defbetinf}
 \betinf:=\liminf_{n\to \infty} \frac{1}{\etnpd}\EE\lt[\Wb^p_Q(\mu_n,n\one_{Q})\rt].
\end{equation}
While it is probably common knowledge in the community that also $\betinf\in (0,\infty)$ (the proof for $d\ge 3$ can be easily adapted from \cite{FoGu15}), we provide a proof in  Proposition \ref{prop:betainf} in the  Appendix.
\begin{remark}
 As for $\betsup$, when $p=d=2$ or $d\ge 3$, the liminf is actually a limit, see \cite{ambrosio2022quadratic,GolTrecombi} or a very simple proof in the case $d\ge 3$ in Proposition \ref{prop:betainf}. Again, this will play no role in our analysis.
\end{remark}
\begin{remark}
 While it is conjectured that $\betinf=\betsup$ for every $p,d\ge 1$ it is currently only known for $p=2$, see \cite{ambrosio2022quadratic,GHT}. See also \cite{goldman2023concave} for the case $d=1$ and $p\in(0,1/2)$.
\end{remark}

For $q>0$ we define the rate function 
\begin{equation}\label{deftaun}
 \tnpd:=\begin{cases}
                \etnpd & \textrm{if } p\in[1,d), \\
                (\log n)^{1+\frac{2}{q}} &\textrm{if } p=d=2,\\
                (\log n)^{\frac{d}{q}} & \textrm{if } d\ge 3 \textrm{ and } p=d,\\
                (\log n)^{d-1-p\lt(1-\frac{1}{q}\rt)} & \textrm{if } d\ge 2 \textrm{ and } p>d.
               \end{cases}
\end{equation}
Let now $V$ be a radial  and $C^1$ function satisfying for some $C>0$ and every $t\ge 1$ (with a slight abuse of notation we identify $V(x)$ and $V(|x|)$),
\begin{equation}\label{hypV}
 \frac{1}{C} t^{q-1}\le V'(t)\le C t^{q-1}.
\end{equation}
We assume moreover that $V(0)$ is chosen such that $\rho(x):=\exp(-V(x))$ is a probability density on $\R^d$. Let $R_n$ be the unique solution to
\begin{equation}\label{defRnmainstat}
 V(R_n)=\log n
\end{equation}
and set 
\begin{multline}\label{defell}
 \ellsup:=\limsup_{n\to \infty} \frac{R_n^d}{(\log n)^{\frac{d}{q}}}\int_{B_1}\lt(1-\one_{d=2} \frac{V(R_n x)}{\log n}\rt)\\ \textrm{and } \qquad \ellinf:=\liminf_{n\to \infty} \frac{R_n^d}{(\log n)^{\frac{d}{q}}}\int_{B_1}\lt(1-\one_{d=2} \frac{V(R_n x)}{\log n}\rt).
\end{multline}
\begin{remark}\label{rem:ell}
 Under hypothesis \eqref{hypV} it is not hard to check that $0<\ellinf\le \ellsup<\infty$. Moreover, it is easy to construct examples of potentials $V$ satisfying \eqref{hypV} but for which $\ellinf<\ellsup$. Let us also notice that the value of $\ellinf, \ellsup$ does not change if instead of \eqref{defRnmainstat} we use
 \[
  V(R_n)+\alpha \log R_n =\log n
 \]
for some $\alpha\in \R$.
\end{remark}
\begin{remark}\label{remtypicalV}
 If $V(t)= V(0)+ t^q/q$ (the case $q=2$ corresponding to the Gaussian) then to leading order,
 \[
  R_n=(q \log n)^{\frac{1}{q}} \qquad \textrm{and } \qquad \ellinf=\ellsup= q^{\frac{d}{q}}\int_{B_1}  (1-\one_{d=2} |x|^q).
 \]
\end{remark}

Our main result is the following:
\begin{theorem}\label{mainthm}
Let $(X_i)_{i=1}^\infty$ be i.i.d. random variables of law $\rho$. Setting $\mu_n:=\sum_{i=1}^n \delta_{X_i}$, we have for $p\le d$,
  \begin{multline}\label{mainstatpled}
\betsup \lt(\int_{\R^d} \rho^{1-\frac{p}{d}} \one_{p<d} +\ellsup \one_{p=d}\rt)\ge \limsup_{n\to \infty} \frac{1}{\tnpd} \EE\lt[W_p^p(\mu_n,n\rho)\rt]\\
\ge \liminf_{n\to \infty}\frac{1}{\tnpd} \EE\lt[W_p^p(\mu_n,n\rho)\rt]\ge \betinf \lt(\int_{\R^d} \rho^{1-\frac{p}{d}} \one_{p<d} +\ellinf \one_{p=d}\rt)
 \end{multline}
and for $p>d$,
\begin{equation}\label{mainstatpegd}
 \EE\lt[W_p^p(\mu_n,n\rho)\rt]\simeq \tnpd.
\end{equation}

\end{theorem}
Theorem \ref{mainthm} is obtained as a combination of Proposition \ref{prop:lowerbound} for the lower bound and Proposition \ref{prop:upperbound} for the upper bound.
\begin{remark}
On the one hand, when $p=d=q=2$ (where $\betsup=\betinf=1/(4\pi)$) we get from Remark \ref{remtypicalV},
\[
 \betinf \ellinf=\betsup\ellsup=\frac{1}{4},
\]
recovering the result from \cite{caglioti2023random}. On the other hand, for $p<d/2$ with $d\ge3$, Theorem \ref{mainthm} is covered by the results in \cite{BaBo,DeScSc13}.
\end{remark}
In Theorem \ref{thm:concentration}, we complement the convergence in expectation with concentration estimates when $p\le d$ and $\rho$ satisfies a Poincar\'e inequality. Similar concentration results were obtained in this range of parameters (actually $p\le 2$ if $d=2$ and  $p<d$ if $d\ge 3$) in the case of the matching to the uniform measure on the cube in \cite[Remark 4.7]{AmStTr16}, see also \cite[Remark 6.5]{goldman2021convergence} or \cite[Proposition 5.3]{GolTrecombi}.

\subsection*{Previous results, heuristics and strategy of proof}
As understood in \cite{BaBo,DeScSc13}, while the matching problem is quite robust when $p<d/2, d\ge 3$, the situation becomes more complicated when $p\ge d/2$ (or $d=2$). While the case of general densities on bounded domains was essentially treated in \cite{BeCa,goldman2021convergence,ambrosio2022quadratic} (see also Theorem \ref{th:upperboundappendix}) the unbounded case remains challenging. Motivated by this question, the study of matching for Gaussian samples was initiated in \cite{Le17}, where using semi-group techniques, Ledoux proved the correct upper bound for the scaling law in the case $p=d=2$ but only a sub-optimal local bound. Soon after the correct lower bound (in terms of scaling laws) was obtained in \cite{talagrand2018scaling}. In that paper, Talagrand proves the correct lower bound $(\log n )^{\frac{d}{2}}$ when $p=d\ge 3$ and suggests that his proofs  should also apply to the family of densities $\rho(x)=\exp(-|x|^q)$ with $q>0$. Still at the level of scaling law, the case $p<d$ with $d\ge 3$ was covered in \cite{ledoux2019optimal}.  Using subadditivity techniques in the spirit of \cite{BeCa,goldman2021convergence,ambrosio2022quadratic}, the optimal constant was recently obtained (still in the Gaussian case) for the case $p=d=2$ in \cite{caglioti2023random}.\\

Here are the main contributions of this paper:
\begin{itemize}
 \item In the Gaussian case we improve on \cite{Le17,talagrand2018scaling,ledoux2019optimal} by obtaining the optimal prefactors. Moreover, we are able to treat (this time only at the level of scaling laws) also the case $p>d$ where no conjecture was even made. A surprising feature of this regime is that the extra logarithmic in dimension $d=2$ disappears, see Remark \ref{rem:logd2}.
 \item We cover a large family of  radial densities thus proving that there is nothing special about the Gaussian measure. In particular \cite{Le17,ledoux2019optimal} strongly relies on the properties of the Ornstein-Uhlenbeck semi-group while \cite{caglioti2023random} used crucially the tensorization properties of the Gaussian measure. Notice that for $q<1$, measures satisfying \eqref{hypV} cannot be log-concave (see also Remark \ref{rem:Poinc}) which shows that also this natural property is not essential. While this goes beyond the scope of this paper, it suggests that for $p<d$, Theorem \ref{mainthm} should hold in the much larger class of densities satisfying $\int_{\R^d} \rho^{1-\frac{p}{d}}<\infty$.
 \item We show that obtaining (sharp) lower bounds in the matching problem is much easier than obtaining the corresponding upper bounds, see Remark \ref{rem:globalp>d}. In particular we include in Appendix \ref{sec:appendix} a short proof of the convergence of the boundary transport problem in the case of the uniform measure on the cube as well as the lower bound on bounded sets. The main new ingredient is the elementary estimate \eqref{eq:change} of Lemma \ref{change}.
 \item For the upper bound we use a general subadditivity property of the Wasserstein distance, see Lemma \ref{subad}. This goes beyond the classical application of subadditivity  in the context of matching problems where usually the measures are required to have disjoint support, see also the discussion below.
\end{itemize}

As in \cite{BeCa,goldman2021convergence,ambrosio2022quadratic,caglioti2023random,GolTrecombi,TreMar}, see also \cite{BaBo,DeScSc13} in the case $p<d/2$, our proof, which we now sketch, relies on subadditivity/superadditivity arguments. We focus on the case $d\ge 3$ to simplify the notation. The idea is to subdivide $\R^d$ in cubes $Q_j=z_j+(0,h_j)^d$ for $h_j>0$ and  $z_j\in \R^d$. Assuming for a moment that for every $j$ (which of course is never the case),
\begin{equation}\label{hypmassabsurd}
\mu_n(Q_j)=n\rho(Q_j),
\end{equation}
we can write by subadditivity (see \eqref{notWOm} for the notation $W_\Omega$),
\[
 W_p^p(\mu_n,n\rho)\le \sum_j W_{Q_j}^p(\mu_n,n\rho).
\]
Now provided $V$ does not oscillate too much on $Q_j$ i.e.
\begin{equation}\label{hyposcillsmall}
 h_j\sup_{Q_j} |\nabla V|\simeq h_j |z_j|^{q-1}\ll1,
\end{equation}
we may consider that in each $Q_j$ the points $X_i$ follow a uniform distribution. Finally, provided there are in expectation many points in each of the cubes i.e.
\begin{equation}\label{hypmanypoints}
 \EE[\mu_n(Q_j)]=n\rho(Q_j)\simeq nh_j^d \rho(z_j)\gg 1,
\end{equation}
we can use the definition of $\betsup$ to infer
\[
 \EE\lt[W_{Q_j}^p(\mu_n,n\rho)\rt]\le (1+o(1))\betsup h_j^p (n\rho(Q_j))^{1-\frac{p}{d}}
 =(1+o(1))\betsup n^{1-\frac{p}{d}} \int_{Q_j} \rho^{1-\frac{p}{d}}.
\]
This estimate is then summed over all cubes satisfying \eqref{hyposcillsmall} and \eqref{hypmanypoints}. Let $R_n:=\max |z_j|$ and set $h=R_n^{1-q}$ so that \eqref{hyposcillsmall} is almost satisfied. From \eqref{hypmanypoints} we see that we must impose
\begin{equation*}\label{hypcondRn}
 n R_n^{d(1-q)} \exp{-V(R_n)}\gg 1.
\end{equation*}
This motivates the choice of  $R_n$ as the unique solution to
\begin{equation}\label{Rnsketch}
 V(R_n)+ d(q-1)\log R_n =\log n.
\end{equation}
Notice that with this definition, we have
\[
 R_n\simeq (\log n)^{\frac{1}{q}}.
\]
Since (see Remark \ref{rem:ell} and the computation \eqref{estim:integrhonp>d})
\begin{multline*}
 n^{1-\frac{p}{d}}\int_{B_{R_n}} \rho^{1-\frac{p}{d}}\\
 =(1+o(n))\tnpd\lt(\one_{p<d} \int_{\R^d} \rho^{1-\frac{p}{d}}+\one_{p= d} \frac{R_n^d}{(\log n)^{\frac{d}{q}}}\int_{B_1}\lt(1-\one_{d=2} \frac{V(R_n x)}{\log n}\rt)+ C\one_{p>d}\rt)
\end{multline*}
we can then conclude the proof of \eqref{mainstatpled} and \eqref{mainstatpegd}.
\begin{remark}\label{rem:logd2}
 In the case $d=2<p$, the integral $\int_{B_{R_n}} \rho^{1-\frac{p}{2}}$ concentrates around $\partial B_{R_n}$ where in each cube $Q_j$ the number of points is of order one. This explains why there is no additional logarithms with respect to the case $d\ge 3$.
\end{remark}
One simple but important observation of this work is that the sketch of proof above can be made rigorous for the lower bound without many changes, see Section \ref{sec:lowerbound}. Indeed, for the boundary functional we do not need to impose condition \eqref{hypmassabsurd} and can locally modify the mass thanks to Lemma \ref{change}. This leads to much simpler proofs than in earlier works on the subject where the boundary  and  standard Wasserstein distances were treated in a similar way, see \cite{ambrosio2022quadratic,GolTrecombi,caglioti2023random}. While this is mostly of esthetical nature in the case $p\le d$, it turns out to be an important observation to treat the case $p>d$, see Remark \ref{rem:globalp>d} below.\\

For the proof of the upper bound, many changes need to be made to the sketch of proof above. First, as already understood in \cite{Le17}, in order to hope for \eqref{hyposcillsmall} and \eqref{hypmanypoints} to hold we need to truncate the measure $\rho$ on the ball $B_{R_n}$. While the construction of \cite{Le17} works for $p<d$ (and also for $p=d=2$ thanks to the extra logarithm, see \cite{caglioti2023random}), it yields a transport cost which is too large for $p\ge d$. As alternative we essentially project the points outside $B_{R_n}$ on $\partial B_{R_n}$, see Lemma \ref{distgamman}. In particular, this truncation imposes the choice
\[
 V(R_n) +\alpha \log R_n =\beta \log n
\]
where (compare with \eqref{Rnsketch})
\begin{equation}\label{defqalpha}
 \beta\begin{cases}\in (\frac{p}{d},1)& \textrm{if } p<d \\ =1 &\textrm{if } p\ge d\end{cases}
 \qquad \textrm{ and } \qquad  \alpha \begin{cases}
                                      =0 &\textrm{if } p<d\\
                                    \in (d(q-1),2q+d(q-1)) &\textrm{if } p=d=2\\
                                              \in (d(q-1),q+d(q-1)) &\textrm{if } p=d\ge3\\
  =d(q-1) &\textrm{if } p>d.
                                                                                         \end{cases}
\end{equation}
The second and most important modification is due to the condition \eqref{hypmassabsurd}. Following \cite{goldman2021convergence} a way around it is to use triangle inequality and Young, see also Lemma \ref{subad} and write for $\eps\in(0,1)$,
\begin{multline}\label{Youngintro}
 \EE[W_p^p(\mu_n,n\rho)]\le (1+ \eps)\sum_{j} \EE\lt[W_{Q_j}^p\lt(\mu_n, \frac{\mu_n(Q_j)}{\rho(Q_j)}\rho\rt)\rt]\\
 + \frac{C}{\eps^{p-1}} \EE\lt[W_p^p\lt(\sum_j \frac{\mu_n(Q_j)}{\rho(Q_j)}\one_{Q_j}\rho, n\rho\rt)\rt].
\end{multline}
The first term, which we sometimes call the local term, may be estimated following now exactly the sketch of proof while the difficulty is to prove that the second term, sometimes called the global term, is negligible. Since the estimate of the global term is usually delicate, see \cite{ambrosio2022quadratic, GolTrecombi,TreMar} and since we need to be very precise in particular for $p\ge d$, we had to resolve to modify the partition in order to take advantage of the radial symmetry of the  measure\footnote{Let us point out that for $q<1$ there is anyway the issue that in order to satisfy \eqref{hyposcillsmall} and \eqref{hypmanypoints} at the same time it is not possible to take $h=h_j$ independent of $Q_j$.}. A crucial requirement for the strategy to work is that all the elements of the partition have the same shape so that the rates of convergence in Theorem \ref{upperboundappendix} are uniform. To this aim we relax the condition of being a partition, see Lemma \ref{subad}, and use instead coverings by sets of the form
\[
 \Delta^{\delta,z}_r=r\{ 1\le |x|\le 1+\delta \ : \hat{x}\in D_\delta(z)\}
\]
where $\hat{x}=x/|x|$ and  $D_\delta(z)$ denotes the image through the exponential map on $\bS^{d-1}$ of a $(d-1)$-dimensional cube on the tangent -- the precise reason for such choice is  technical, since it is common in the literature to argue with cubes. The estimate of the corresponding global term in Proposition \ref{prop:global} is divided into a radial, see Section \ref{sec:radial} and an angular, see Section \ref{sec:angular}, transport. While for most values of $p$ and $d$ the order of subdivision (first radial and then angular or vice-versa) does not matter, for certain values it seems instead crucial (see Remark \ref{rem:diff}).
\begin{remark}
 Let us point out that a similar strategy for the upper bound was suggested in \cite{talagrand2018scaling} for the case $p=d=2$.
\end{remark}

\begin{remark}\label{rem:globalp>d}
 In the case of the lower bound, one could argue as in the upper bound and use triangle inequality together with Young to obtain the analog of \eqref{Youngintro}. The main difference would be the minus sign in front of the global term. Using that $\Wb\le W$, one could then use the estimates for the upper bounds to conclude. This would however not work when $p>d$ since the local and the global terms are both of the order of $\tnpd$, see Proposition \ref{prop:global}.
\end{remark}

\section{Notation and preliminary results}
\subsection{Notation}
We identify measures which are absolutely continuous with respect to the Lebesgue measure with their density. Moreover, when integrating against the Lebesgue measure, and when  it is clear from the context, we omit both integration variables and the measure in the integrals. We write $|\Omega|$ for the Lebesgue measure of a set $\Omega$. We write $\HH^k$ for the Hausdorff measure of dimension $k$.
We denote by $\mathcal{M}(\Omega)$ the set of signed measures on $\Omega$ and by $\mathcal{M}^+(\Omega)$ the set of the positive measures on $\Omega$.\\
The notation $A\les B$, which we use in output statements, means that there exists a universal (meaning depending possibly on $d$ and $p$)
constant $C>0$ such that $A\le C B$. We write $A\les_\eps B$ if the implicit constant depends on the parameter $\eps$. We write $A\simeq B$ if $A\les B$ and $B\les A$. We denote by  $\omega(t)$ a rate function which may vary from line to line such that $\lim_{t\to 0} \omega(t)=\lim_{t\to \infty} \omega(t)=0$ We write $\omega_\eps$ to indicate that the rate function also depends on $\eps$.\\
For $d\geq2$ we write  $\bS^{d-1}=\partial B_1$ for the unit sphere and denote by
$$ \sigma:= \frac{1}{\HH^{d-1}(\bS^{d-1})} \HH^{d-1}\restr \bS^{d-1}$$
the uniform probability measure on $\mathbb{S}^{d-1}$.
Moreover $\mathbb{S}^{d-1}$ is endowed with the (Riemannian) spherical distance $\mathsf{d}_{\mathbb{S}^{d-1}}$. We write $\bStilde^{d-1} := \bS^{d-1}\times O(d-1)$, where $O(d-1)$ denotes the orthogonal group on $\R^{d-1}$. For $z = (y, U) \in \bStilde^{d-1}$ where $p \in \mathbb{S}^{d-1}$, $U \in O(d-1)$, we define
\begin{equation}\label{defDdelta}
 D_{\delta}(z):=\exp_{y}( U (-\delta/2,\delta/2)^{d-1}) \subseteq \bS^{d-1}.
 \end{equation}
 where $\exp_{y}$ denotes the exponential map at $y$. For $z = (e_1, \operatorname{Id}_{d-1})$ we simply write $D_{\delta}$. Notice that we are implicitly identifying the tangent at $y$ with $\R^{d-1}$, which can is done e.g.\ by viewing $\bS^{d-1}$ as a submanifold of $\R^{d-1}$. We also write $\stilde$ for the uniform (Haar) probability measure on $\bStilde$. 

\subsection{The Wasserstein distance}
For  $\mu,\lambda\in \mathcal{M}^+(\R^d)$ with $\mu(\R^d)=\lambda(\R^d)$, we define the $p-$Wasserstein distance  between $\mu$ and $\lambda$ as
\begin{eqnarray}\label{wass}
W_{p}^p(\mu,\lambda):=\left.\min\left\{\int_{\R^d\times\R^d}|x-y|^pd\pi(x,y)\right|\pi_1=\mu,\pi_2=\lambda\right\}.
\end{eqnarray}
For a set $\Omega$ we write 
\begin{equation}\label{notWOm}
 W_{\Omega}^p(\mu,\lambda):=W_{p}^p(\mu\restr \Omega,\lambda\restr \Omega).
\end{equation}
We recall the following classical consequence of the Benamou-Brenier formula, see e.g. \cite{peyre2018comparison} or \cite[Lemma 3.4]{goldman2021convergence}.

\begin{proposition}\label{BBcomodo} Let $\Omega$ be an open set and $\mu,\lambda\in \mathcal{M}^+(\Omega)$ be such that  $\mu(\Omega)=\lambda(\Omega)$. If $\phi$ is a solution of 
\[
\Delta\phi=\mu-\lambda \quad \textrm{in } \Omega \qquad \textrm{and } \qquad \nabla\phi\cdot\nu=0 \quad  \textrm{on } \partial\Omega,
\]
where $\nu$ denotes the exterior normal to $\partial\Omega$, then
\begin{eqnarray*}
W_p^p(\mu,\lambda)\lesssim\int_{\Omega}\frac{|\nabla\phi|^p}{\mu^{p-1}}.
\end{eqnarray*}
\end{proposition}
We now state a simple subadditivity result for the Wasserstein distance whose proof is immediate. It extends more classical subadditivity results to the case of coverings which are not necessarily partitions.   The last statement of the Lemma comes from combining the second one with the property $(a+b)^p\leq(1+\varepsilon)a^p+\frac{c}{\varepsilon^{p-1}}b^p$.
\begin{lemma}\label{subad}
Let $(E,{\mathcal E},\zeta)$ be a measure space and $(\mu_z,\lambda_z)_{z\in E}$ a collection of measures on $\Omega$ such that $\mu_z(\Omega)=\lambda_z(\Omega)$. Then
$$
W_p^p\left(\int_E\mu_zd\zeta(z),\int_E\lambda_zd\zeta(z)\right)\leq\int_EW_p^p(\mu_z,\lambda_z)d\zeta(z).
$$
In particular, for $\Omega\subseteq\mathbb{R}^d$, if  $(\Omega_z)_{z\in E}$ is a collection of Borel subsets of $\Omega$,  such that
\begin{eqnarray*}
\int_E\mathbbm{1}_{\Omega_z}d\zeta(z)=1\qquad\lambda+\mu-a.e. \qquad \textrm{and } \qquad \mu(\Omega_z)=\lambda(\Omega_z) \textrm{ for } \zeta\ a.e. z,
\end{eqnarray*}
then
\begin{eqnarray*}
W_p^p(\mu,\lambda)\leq\int_EW_{\Omega_z}^p(\mu,\lambda)d\zeta(z).
\end{eqnarray*}
Therefore, letting $\kappa:= \mu(\Omega)/\lambda(\Omega)$ and $\kappa_z:= \mu(\Omega_z)/\lambda(\Omega_z)$, there exists $c=c(p,d)>0$ such that for every $\varepsilon\in(0,1)$
\begin{equation}\label{locglob}
 W_p^p\left(\mu,\kappa\lambda\right)\leq(1+\varepsilon)\int_EW_{\Omega}^p\left(\mu,\kappa_z \lambda\right)d\zeta(z)
+\frac{c}{\varepsilon^{p-1}}W_p^p\left(\int_E\kappa_z\lambda\mathbbm{1}_{\Omega_z}d\zeta(z),\kappa\lambda\right).
\end{equation}
\begin{remark}
To recover classical applications of subadditivity, see e.g. \cite[Lemma 3.1]{goldman2021convergence}, it is enough to apply the lemma  with $E$ finite and $\zeta=\sum_{z} \delta_z$ and obtain 
$$
W_p^p\left(\sum_{z\in E}\mu_z,\sum_{z\in E}\lambda_z\right)\leq\sum_{z\in E}W_p^p(\mu_z,\lambda_z).
$$
\end{remark}
\end{lemma}
\subsection{The boundary Wasserstein distance}
For the lower bound we will use the boundary Wasserstein distance, that has been defined in \cite{FG2010} between $\mu$ and $\lambda$ measures on $\Omega$ not necessarily of the same mass. It is defined as 
\begin{eqnarray}\label{bwass}
\Wb_{\Omega}^p(\mu,\lambda):=\left.\min\left\{\int_{\bar\Omega\times\bar\Omega}|x-y|^pd\pi(x,y)\right|\pi_1\lfloor\Omega=\mu,\pi_2\lfloor\Omega=\lambda\right\}.
\end{eqnarray}
As for $W_\Omega$ we write for $\Omega'\subset \Omega$, 
\[
 \Wb_{\Omega'}^p(\mu,\lambda):=\Wb_{\Omega'}^p(\mu\restr \Omega',\lambda\restr\Omega').
\]

We recall that $\Wb$ is superadditive, see e.g. \cite{ambrosio2022quadratic} or \cite[Lemma 25]{BaBo}.

\begin{lemma}\label{superad}
Let $\Omega\subset \R^d$ be an open set and let $\Omega_j\subset \Omega$ be pairwise disjoint open sets then
\begin{eqnarray}\label{superad2}
\Wb_\Omega^p(\mu,\lambda)\geq \sum_j \Wb^p_{\Omega_j}(\mu,\lambda).
\end{eqnarray}
\end{lemma}

We conclude with a simple lemma which will allow us to adjust the densities when working with $\Wb$.
\begin{lemma}\label{change}
Let $\Omega$ be a bounded set. Then for every $m,m'>0$
\begin{equation}\label{eq:change}
\Wb_\Omega^p\lt(m,m'\rt)\lesssim \diam(\Omega)^{p+d}  \frac{|m-m'|^p}{\max(m,m')^{p-1}}.
\end{equation}
\end{lemma}
\begin{proof} 
Since $\Wb_{\Omega}\le \Wb_{{\rm Conv}\,  \Omega}$ where ${\rm Conv}\,  \Omega$ denotes the convex enveloppe of $\Omega$, we may assume without loss of generality that $\Omega$ is convex. Let us first treat the one dimensonal case. We may assume that $\Omega=(0,L)$ for some $L>0$ and $m\ge m'$. We then set for $x\in (0,L)$,
\[
 T(x):=\min\lt(\frac{m}{m'}x ,L\rt)
\]
so that $T$ sends $m \one_{(0,L)}$ on $m'\one_{(0,L)}+ (m-m')/L \delta_{L}$ and is thus $({\rm Id}\times T)_\sharp (m\one_{(0,L)})$ is an admissible transport plan for $\Wb_{\Omega}^p(m,m')$. We then compute
\[
 \Wb_{\Omega}^p(m,m')\le m\int_0^L |T-x|^p =m\int_0^{\frac{m'}{m} L} \lt|\frac{m}{m'}-1\rt|^p x^p + m\int_{\frac{m'}{m} L}^L (L-x)^p.
\]
An elementary computation then yields the claimed
\begin{equation}\label{Wb1d}
 \Wb_{\Omega}^p(m,m')\les L^{p+1} \frac{|m-m'|^p}{m^{p-1}}.
\end{equation}
Let us now conclude the proof in the case $d\ge 2$. For $x'\in \R^{d-1}$ let
\[
 \Omega_{x'}=\{t\in \R : \ (t,x')\in \Omega\}
\]
be the intersection of $\Omega$ with the line $\R e_1+ x'$. Since $\Omega$ is convex this is a segment. Letting $\Omega':=\{x'\in \R^{d-1} \, : \, \Omega_{x'}\neq \emptyset\}$ and for $x'\in \Omega_{x'}$, $\pi_{x'}$ be an admissible  transport plan for $\Wb_{\Omega_{x'}}^p(m,m')$ we see that the plan defined by
\[
 \int_{\overline{\Omega}\times \overline{\Omega}} \zeta d\pi:=\int_{\Omega'} \int_{\overline{\Omega}_{x'}\times \overline{\Omega}_{x'}} \zeta(x'+te_1,x'+se_1)d\pi_{x'}(t,s) dx'
\]
is admissible for $\Wb_{\Omega}^p(m,m')$ and thus
\begin{equation*}\begin{split}
 \Wb_{\Omega}^p(m,m') & \le \int_{\Omega'} \Wb_{\Omega_{x'}}^p(m,m') dx'\stackrel{\eqref{Wb1d}}{\les} \int_{\Omega'} \diam(\Omega_{x'})^{p+1} \frac{|m-m'|^p}{m^{p-1}}\\
 & \les \diam(\Omega)^{p+d} \frac{|m-m'|^p}{m^{p-1}}.
 \end{split}
\end{equation*}
This concludes the proof of \eqref{eq:change}.
\end{proof}

\section{Proof of the lower bound}\label{sec:lowerbound}
In this section we prove the lower bound part of Theorem \ref{mainthm}. Let us recall the definition \eqref{deftaun} of $\tnpd$ and \eqref{defell} of $\ellinf$.

\begin{proposition}\label{prop:lowerbound}
 Let $(X_i)_{i=1}^\infty$ be i.i.d. random variables of law $\rho$. Setting $\mu_n:=\sum_{i=1}^n \delta_{X_i}$, we have for $p\le d$,
 \begin{equation}\label{mainstatpledlower}
 \liminf_{n\to \infty}\frac{1}{\tnpd} \EE\lt[W_p^p(\mu_n,n\rho)\rt]\ge \betinf \lt(\int_{\R^d} \rho^{1-\frac{p}{d}} \one_{p<d} + \ellinf \one_{p=d}\rt)
 \end{equation}
and for $p>d$,
\begin{equation}\label{mainstatpegdlower}
 \EE\lt[W_p^p(\mu_n,n\rho)\rt]\ges \tnpd.
\end{equation}

\end{proposition}
\begin{proof}
{\it Step 1.}
Let us first set some notation. For $h>0$ and  $z\in h\Z^d$, we denote $Q_z:= z+(0,h)^d$. For any $Z\subset h\Z^d$, we set
\begin{equation}\label{defOmegaZ}
 \Omega_Z:=\cup_{z\in Z} Q_z.
\end{equation}
Let $\eps>0$ be fixed. We first claim that provided there exists $m=m_\eps$ large enough such that for $z\in Z$,
\begin{equation}\label{goodbounds}
h V'(|z|)\les \eps
\end{equation}
and
\begin{equation}\label{averagepoints}
   n h^d \rho(|z|)\ges_\eps \exp(C \eps m)
\end{equation}
then
\begin{equation}\label{claimsuperad}\begin{split}
 & \EE\lt[W^p_p(\mu_n,n\rho)\rt]\\
 & \ge(1-\omega(\eps)-\max_{z\in Z}\omega(n\rho(Q_z)))\betinf n^{1-\frac{p}{d}}\sum_{z\in Z}\int_{Q_z} \rho^{1-\frac{p}{d}}(1+\one_{d=2} (\log n \rho(Q_z))^{\frac{p}{2}}).
 \end{split}
\end{equation}
Let us first prove the claim. To this aim we first use superadditivity, see Lemma \ref{superad}, to obtain
\[
\EE\lt[W^p_p(\mu_n,n\rho)\rt]\ge \sum_{z\in Z} \EE\lt[\Wb_{Q_z}^p(\mu_n,n\rho)\rt].
\]
Fix $z\in Z$ and set $N:=\mu_n(Q_z)$. Let $\hat{X}_i$ be i.i.d. random variables of law $\hat{\lambda}:=\rho\one_{Q_z}/\rho(Q_z)$ so that $\mu_n\restr Q_z$ has the same law as $\hat{\mu}_N:=\sum_{i=1}^N \delta_{\hat{X}_i}$. We thus have
\[
 \EE\lt[\Wb_{Q_z}^p(\mu_n,n\rho)\rt]=\EE\lt[\EE\lt[\Wb_{Q_z}^p(\hat{\mu}_k,n\rho)|N=k\rt]\rt].
\]
By \eqref{goodbounds}, we can use either the Knothe map as in \cite[Lemma 1]{BeCa} or \cite[Proposition 2.4]{ambrosio2022quadratic} (in both cases the proof easily extend to general dimensions) to find a Lipschitz map $T: Q_z\to Q_z$ such that $T_\sharp \hat{\lambda}=1/|Q_z|$ and such that $\Lip T\le 1+\omega(\eps)$. Letting $\mu'_k:=\sum_{i=1}^k \delta_{T(\hat{X}_i)}$ we find for every $k$,
\[
 \EE\lt[\Wb_{Q_z}^p(\hat{\mu}_k,n\rho)\rt]\ge (1-\omega(\eps))\EE\lt[\Wb_{Q_z}^p\lt(\mu'_k,\frac{n\rho(Q_z)}{|Q_z|}\rt)\rt].
\]
Using triangle inequality and \eqref{eq:change} from Lemma \ref{change}, we find
\begin{multline*}
\Wb_{Q_z}^p(\hat{\mu}_k,n\rho)\ge (1-\omega(\eps))\Wb_{Q_z}^p\lt(\mu'_k,\frac{k}{|Q_z|}\rt) -\frac{C}{\eps^{p-1}}\Wb_{Q_z}^p\lt(\frac{k}{|Q_z|},\frac{n\rho(Q_z)}{|Q_z|}\rt)\\
\ge (1-\omega(\eps))\Wb_{Q_z}^p\lt(\mu'_k,\frac{k}{|Q_z|}\rt) -\frac{C}{\eps^{p-1}} |Q_z|^{\frac{p}{d}} \frac{|k-n\rho(Q_z)|^p}{(n\rho(Q_z))^{p-1}}.
\end{multline*}
Taking expectation and using the definition \eqref{defbetinf} of $\betinf$ (recall also the definition \eqref{defetanpd} of $\etnpd$ we find
\[
\EE\lt[\Wb_{Q_z}^p(\hat{\mu}_k,n\rho)\rt]\ge (1-\omega(\eps)-\omega(k))\betinf |Q_z|^{\frac{p}{d}} \eta_k-\frac{C}{\eps^{p-1}} |Q_z|^{\frac{p}{d}} \frac{|k-n\rho(Q_z)|^p}{(n\rho(Q_z))^{p-1}}.
\]
Taking the expectation over $N$
and using that $\EE[|N-n\rho(Q_z)|^p]\les n\rho(Q_z)^{p/2}$ (recall that \eqref{averagepoints} holds and thus $\EE[N]=n\rho(Q_z)\gg1$ provided $m$ is large enough) we get
\begin{multline*}
 \EE\lt[\Wb_{Q_z}^p(\mu_n,n\rho)\rt]\ge (1-\omega(\eps)-\omega(n\rho(Q_z)))\betinf |Q_z|^{\frac{p}{d}} \EE[\eta_N]-\frac{C}{\eps^{p-1}} |Q_z|^{\frac{p}{d}} (n\rho(Q_z))^{1-\frac{p}{2}}\\
 = (1-\omega(\eps)-\omega(n\rho(Q_z)))\betinf |Q_z|^{\frac{p}{d}} \eta_{n\rho(Q_z)}-\frac{C}{\eps^{p-1}} |Q_z|^{\frac{p}{d}} (n\rho(Q_z))^{1-\frac{p}{2}}.
\end{multline*}
Recalling that
\[
 \eta_{n\rho(Q_z)}=(n\rho(Q_z))^{1-\frac{p}{d}}(1+\one_{d=2} (\log n \rho(Q_z))^{\frac{p}{2}} )
\]
we see that provided $m$ is large enough (depending on $\eps$) we have
\[
 \frac{1}{\eps^{p-1}}(n\rho(Q_z))^{1-\frac{p}{2}}\ll \eta_{n\rho(Q_z)}
\]
so that
\begin{multline*}
  \EE\lt[\Wb_{Q_z}^p(\mu_n,n\rho)\rt]\ge (1-\omega(\eps)-\omega(n\rho(Q_z)))\betinf |Q_z|^{\frac{p}{d}} \eta_{n\rho(Q_z)}\\
  \ge (1-\omega(\eps)-\omega(n\rho(Q_z)))\betinf n^{1-\frac{p}{d}}\int_{Q_z} \rho^{1-\frac{p}{d}}(1+\one_{d=2} (\log (n \rho(Q_z)))^{\frac{p}{2}} ).
\end{multline*}
After summation this proves \eqref{claimsuperad}.\\

{\it Step 2.}
We now consider separately the cases $p<d$, $p=d$ and $p>d$.\\
{\it Step 2.1.} The case $p<d$. Let $R>0$ be fixed, set
\[h=\eps\min(1, R^{1-q}) \qquad \textrm{and} \qquad   Z:=\{z\in h\R^d \ : |z|\le R\}.\]
In this case \eqref{goodbounds} holds and
\[
 \min_Z(n h^d \rho(|z|))\ges_R n \qquad \textrm{and}  \qquad \log (n \rho(Q_z))\ge \log n- C_R
\]
so that \eqref{averagepoints} holds. Applying \eqref{claimsuperad} we find (recall the definition \eqref{deftaun} of $\tnpd$)
\[
 \EE\lt[W^p_p(\mu_n,n\rho)\rt]\ge(1-\omega(\eps)-\omega_{\eps,R}(n))\betinf \tnpd \int_{\Omega_Z} \rho^{1-\frac{p}{d}}.
\]
Dividing by $\tnpd$ and sending $n\to \infty$ we find
\[
 \liminf_{n\to \infty} \frac{1}{\tnpd} \EE\lt[W^p_p(\mu_n,n\rho)\rt]\ge (1-\omega(\eps))\int_{B_R}\rho^{1-\frac{p}{d}}.
\]
Sending $\eps\to 0$ and $R\to \infty$ concludes the proof of \eqref{mainstatpledlower} in this case.\\

{\it Step 2.2.} In the case $p=d$ we choose $\alpha>\max(d(q-1),0)$ and set $R_n$ to be the unique solution to
\[
V(R_n)+\alpha \log R_n=\log n.
\]
We then set
\[
 h:=\eps\min(1, R_n^{1-q}) \qquad \textrm{and} \qquad Z:=\{z\in h\Z^d \ : \ Q_z\cap B_{R_n}\neq\emptyset\}.
\]
By definition of $h$ and \eqref{hypV}, we see that \eqref{goodbounds} holds. Let us check that also \eqref{averagepoints} holds. For this we see that for $z\in Z$,
\[
 n h^d\rho(|z|)\ges n h^d \rho(R_n)=\eps n \min(1, R_n^{d(1-q)}) \exp(-V(R_n))= \eps  \min(1, R_n^{d(1-q)}) R_n^\alpha
\]
so that by our choice of $\alpha$ we get that \eqref{averagepoints} holds and
so that
\[
 \max_{z\in Z} \omega(n\rho(Q_z))=\omega_\eps(n).
\]
By \eqref{claimsuperad} and the fact that for $d=2$, 
\[
 \sum_{z\in Z} \int_{Q_z} |\log |Q_z||\les |\Omega_Z| \log R_n \les |B_{R_n}|\log R_n\ll \tnpd
\]
we have
\begin{multline*}
  \EE\lt[W^p_p(\mu_n,n\rho)\rt]\ge(1-\omega(\eps)-\omega_{\eps}(n))\betinf \int_{\Omega_Z} (1+\one_{d=2} \log(n \rho))\\
  \ge (1-\omega(\eps)-\omega_{\eps}(n))\betinf \int_{B_{R_n}} (1+\one_{d=2} \log(n \rho)).
\end{multline*}
If $d\ge 3$, since $\tnpd=(\log n)^{\frac{d}{q}}$ this concludes the proof of \eqref{mainstatpledlower}. If instead $d=2$, we write 
\[
 \int_{B_{R_n}} \log(n \rho)= \int_{B_{R_n}} (\log n -V(x))=   R_n^d \log n \int_{B_1} \lt(1 -\frac{V(R_n x)}{\log n}\rt)
 \]
 concluding similarly the proof of \eqref{mainstatpledlower}.
\\
{\it Step 2.3.} The case $p>d$. We finally choose $R_n$ to be the unique solution to
\begin{equation}\label{defVnlowerbound}
 V(R_n)+ d(q-1) \log R_n=\log n
\end{equation}
and set for $m\gg1$,
\[
 h:=\eps R_n^{1-q} \qquad \textrm{and } \qquad Z:=\{z\in h\Z^d \ : \ |z|\in [R_n -4 m h,R_n-mh]\}.
\]
Since for $z\in Z$, $|z|\simeq R_n$ as above \eqref{goodbounds} is easily seen to hold. Let us check that also \eqref{averagepoints} is satisfied. For $z\in Z$ we have
\begin{multline*}
 n h^d \rho(|z|)\ge n h^d \rho(R_n-mh)=n h^d \exp(-V(R_n-mh))\\
 = nh^d \exp(-V(R_n))\exp(\int_{R_n-mh}^{R_n} V'(t))\simeq_\eps \exp(\int_{R_n-mh}^{R_n} V'(t))\\
 \ge \exp(C (R_n^q-(R_n-mh)^q))\ge \exp(Cm\eps).
 \end{multline*}
Thus \eqref{averagepoints} holds. Applying \eqref{claimsuperad} yields
\begin{multline*}
  \EE\lt[W^p_p(\mu_n,n\rho)\rt]\ges_\eps n^{1-\frac{p}{d}}\int_{\Omega_Z} \rho^{1-\frac{p}{d}}\\
  \ge n^{1-\frac{p}{d}}\int_{R_n-3mh}^{R_n-2mh} r^{d-1} \rho^{1-\frac{p}{d}}(r)\ges_m n^{1-\frac{p}{d}} R_n^{d-q} \rho(R_n)^{1-\frac{p}{d}}\\
  =n^{1-\frac{p}{d}} R_n^{d-q}\exp\lt(\lt(\frac{p}{d}-1\rt)V(R_n)\rt)\stackrel{\eqref{defVnlowerbound}}{=}R_n^{d-q -(q-1)(p-d)}\simeq \tnpd.
  \end{multline*}
  This proves \eqref{mainstatpegdlower}.
\end{proof}

\section{Proof of the upper bound}
\subsection{The geometric setup}\label{Sec:geomsetup}
Motivated by  the heuristics from the introduction, we set $r_1\ge 1$ be arbitrarily chosen such that in $[r_1,\infty)$, the function $t^{d-q}\exp(-V(t))$ is decreasing and then for $j\ge 1$, set
\begin{equation}\label{defrj}
   r_{j}:=(1+\delta_j)r_{j-1}
\end{equation}
 where for   $\eps \in (0,1)$,
 \begin{equation}\label{defdelta}
  \delta_j:=\eps\begin{cases}
   r_{j-1}^{-q} & \textrm{if $d=2$ and $p\in [1,2)$}\\
            (\log n)^{-1} &\textrm{otherwise}.
           \end{cases}
 \end{equation}
 For the distinction between the cases see Remark \ref{rem:diff} below.
   It is not hard to check that in both cases $r_j\to \infty$ as $j\to \infty$. When $\delta_j= \eps (\log n)^{-1}$ we simply write $\delta$ for $\delta_j$. Let $\bar R_n$ be the unique solution to 
   \begin{equation}\label{defbarRn}
    V(\bar R_n)+\alpha \log \bar R_n=\beta \log n
   \end{equation}
where we recall that $\alpha$ and $\beta$ are defined in \eqref{defqalpha}. Notice that  by \eqref{hypV}, we have
\[
 \bar R_n\simeq (\log n)^{\frac{1}{q}}.
\]
   Let $j_n$ be the first index such that
\[
 r_{j_n}\ge \bar R_n.
\]
We then set
\begin{equation}\label{defRn}
 R_n:=r_{j_n}, \qquad \Omega_n:=B_{R_n}, \qquad R'_n:=r_{j_n-1}, \quad \textrm{and } \quad \Omega'_n:=B_{R'_n}.
\end{equation}
It is an elementary computation to check that
\begin{equation}\label{condRn}
 |R_n-\bar R_n|\les \eps \bar R_n^{1-q} \qquad \textrm{ and } \qquad \delta_{j_n}\simeq \eps (\log n)^{-1}\simeq \eps R_n^{-q}.
\end{equation}
We define the truncated probability density $\rho_n$ by
\begin{equation}\label{defgamman}
\rho_n:=\lt( \one_{\Omega'_n}+ \frac{1-\rho(\Omega'_n)}{\rho(\Omega_n\backslash \Omega'_n)}\one_{\Omega_n\backslash \Omega'_n}\rt)\rho.
\end{equation}
Let us observe that since $(R_n-R'_n)R_n^{q-1}=\delta_{j_n} R_n^q \stackrel{\eqref{condRn}}{\simeq}\eps$, \eqref{hyposcillsmall} is satisfied and thus
\begin{equation}\label{densexternalannulus}
\frac{1-\rho(\Omega_n')}{\rho(\Omega_n\setminus\Omega_n')}=\frac{\int_{R_n'}^{\infty}r^{d-1}\rho(r)}{\int_{R_n'}^{R_n}r^{d-1}\rho(r)}\simeq \frac{R_n^{d-q}\rho(R_n)}{R_n^{d-1}\rho(R_n)R_n\delta_{j_n}}=\frac{1}{R_n^q\delta_{j_n}}\simeq\frac{1}{\varepsilon}.
\end{equation}
For $0\le r<R\le \infty$, $z\in \bStilde^{d-1}$ and $\delta>0$, we define the set (recall the definition \eqref{defDdelta} of $D_\delta(z)$)
\begin{equation}\label{defLambda}
 \Lambda_{r,R}^{z,\delta}=\{ r\le |x|\le R \ : \hat{x} \in D_\delta(z)\}
\end{equation}
where $ \hat{x}=x/|x|$. When $\delta \ge 2\pi$  so that $D_{\delta}(z) = \bS^{d-1}$, we simply write
\begin{equation}\label{defannulus}
 A_{r,R}:=\{ r\le |x|\le R\}= \Lambda_{r,R}^{z,\delta}.
\end{equation}
Finally for $z\in \bStilde^{d-1}$ and $2\le j\le j_n$, we define the quasi-cylindrical  domains
\begin{equation}\label{defDeltazj}
 \Dzj=\Lambda^{z,\delta_j}_{r_{j-1},r_{j}}=  \{  x \in \R^d\, : \, r_{j-1} \le |x|\le (1+\delta_j)r_{j-1} \quad \text{and} \quad \hat{x} \in D_{\delta_j}(z)\}.
\end{equation}

\subsection{The upper bound}
In this section we prove the upper bound part of Theorem \ref{mainthm}. We recall the definition \eqref{deftaun} of $\tnpd$ and \eqref{defell} of $\ellsup$.
\begin{proposition}\label{prop:upperbound}
 Let $(X_i)_{i=1}^\infty$ be i.i.d. random variables of law $\rho$. Setting $\mu_n:=\sum_{i=1}^n \delta_{X_i}$, we have for $p\le d$,
 \begin{equation}\label{mainstatpledupper}
 \limsup_{n\to \infty}\frac{1}{\tnpd} \EE\lt[W_p^p(\mu_n,n\rho)\rt]\le \betsup \lt(\int_{\R^d} \rho^{1-\frac{p}{d}} \one_{p<d} +\ellsup \one_{p=d}\rt)
 \end{equation}
and for $p>d$,
\begin{equation}\label{mainstatpegdupper}
 \EE\lt[W_p^p(\mu_n,n\rho)\rt]\les \tnpd.
\end{equation}

\end{proposition}
\begin{proof}

Let $T$ be the optimal transport map for $W_p(\rho,\rho_n)$. Let $\hat{X}_i:=T(X_i)$ and then $\hat{\mu}_n=T_\sharp \mu_n=\sum_{i=1}^n \delta_{\hat{X}_i}$. Notice that $\hat{X}_i$ are i.i.d. random variables with law $\rho_n$ and $\hat{\mu}_n\restr \Omega'_n=\mu_n\restr \Omega_n$ (recall \eqref{defRn}). Since
\[
 \EE[W_p^p(\mu_n,\hat{\mu}_n)]\le n W_p^p(\rho,\rho_n),
\]
we have by triangle inequality and Young that for $\eps\in (0,1)$,
\[
 \EE[W_p(\mu_n, n\rho)]\le (1+\eps)\EE[W_p^p(\hat{\mu}_n, n\rho_n)]+ \frac{C}{\eps^{p-1}} n W_p^p(\rho,\rho_n).
\]
Since $n W_p^p(\rho,\rho_n)=o(\tnpd)+ \one_{p>d} O(\tnpd)$ by Lemma \ref{distgamman}, we are only left with estimating $\EE[W_p^p(\hat{\mu}_n, n\rho_n)]$. To this aim,
we set 
\begin{equation}\label{defkappajz}\kappa_{1}:= \frac{\hat{\mu}_n(B_{r_1})}{\rho_n(B_{r_1})} \qquad \textrm{and } \qquad \kappa_j^z:= \frac{\hat{\mu}_n(\Delta_j^z)}{\rho_n(\Delta_j^z)}.\end{equation}
Using the subadditivity inequality \eqref{locglob} from Lemma \ref{subad} on $E = \mathbb{N}\times \bStilde^{d-1}$ with
\[
\zeta:=\sum_{j=1}^{j_n} \frac{ \one_j\otimes d \stilde (z) }{ \stilde_j }, \quad \text{where} \quad \stilde_j = \int_{\bStilde^{d-1}} \one_{D_{\delta_j}(z) }(e_1) d \stilde(z)\]
and
\[ \Omega_{1,z}:=B_{r_1}, \quad \textrm{ and } \quad  \Omega_{j,z}:=\Delta_j^z \textrm{ for  } j\ge 2
\]
 we find
\begin{multline*}
 \EE[W_p^p(\hat{\mu}_n, n\rho_n)]\le (1+\eps)\lt( \EE[W_{B_{r_1}}^p(\hat{\mu}_n, \kappa_1 \rho_n)] +\sum_{j=2}^{j_n} \frac{1}{\stilde_j} \int_{\bStilde^{d-1}} \EE[W_{\Delta_j^z}^p(\hat{\mu}_n, \kappa_j^z \rho_n)] d\stilde(z)\rt)\\
 +\frac{C}{\eps^{p-1}}\EE\lt[W_p^p\lt(\kappa_1 \one_{B_{r_1}}\rho_n+\sum_{j=2}^{j_n}\frac{1}{\stilde_j}\int_{\bStilde^{d-1}}\kappa_j^z \one_{\Delta_j^z} \rho_n d\stilde(z), n\rho_n\rt)\rt].
\end{multline*}
Let us first consider the case $p\le d$. In this case, applying \eqref{eq:proplocaltermssmall} from Proposition \ref{prop:localterms} and \eqref{eq:global} from Proposition \ref{prop:global} we obtain (recall the definition \eqref{defetanpd} of $\etnpd$
\[
 \EE[W_p^p(\hat{\mu}_n, n\rho_n)]\le \lt(1+\omega(\eps)+\omega(n)\rt) \betsup \etnpd\int_{\Omega_n} \rho_n^{1-\frac{p}{d}} \lt(1+\one_{d=2} \lt(\frac{\log n \rho_n}{\log n}\rt)^{\frac{p}{2}}\rt)+o_\eps(\tnpd).
\]
Using that for $p\le d$,
\[
 \etnpd\int_{\Omega_n}  \rho_n^{1-\frac{p}{d}} \lt(1+\one_{d=2} \lt(\frac{\log n \rho_n}{\log n}\rt)^{\frac{p}{2}}\rt)\le(1+\omega_\eps(n))\tnpd \int_{\R^d} \rho^{1-\frac{p}{d}}(\one_{p<d}+ \ellsup \one_{p=d})
\]
this concludes the proof of \eqref{mainstatpledupper} sending first $n\to\infty$ and then $\eps\to 0$.\\
In the case $p>d$, we use \eqref{eq:proplocaltermslarge} from Proposition \ref{prop:localterms} and \eqref{eq:global} from Proposition \ref{prop:global} to obtain 
\[
 \EE[W_p^p(\hat{\mu}_n, n\rho_n)]\les_\eps  n^{1-\frac{p}{d}}\int_{\Omega_n} \rho_n^{1-\frac{p}{d}} +\tnpd.
\]
Since
\begin{multline}\label{estim:integrhonp>d}
 n^{1-\frac{p}{d}}\int_{\Omega_n} \rho_n^{1-\frac{p}{d}}\simeq_\eps n^{1-\frac{p}{d}}\int_0^{R_n} r^{d-1}\rho(r)^{1-\frac{p}{d}}\\
 =n^{1-\frac{p}{d}} \rho(R_n)^{1-\frac{p}{d}} \int_{0}^{R_n} r^{d-1} \exp\lt(-\lt(\frac{p}{d}-1\rt)\int_r^{R_n} V'(t)dt\rt) dr\\
 \stackrel{\eqref{hypV}}{\simeq} (n\rho(R_n))^{1-\frac{p}{d}} \int_{0}^{R_n} r^{d-1} \exp\lt(-C (R_n^q-r^q)\rt) dr\simeq
 R_n^{d-q}(n\rho(R_n))^{1-\frac{p}{d}}\stackrel{\eqref{condRn}}{\simeq} \tnpd,
\end{multline}
this concludes the proof of \eqref{mainstatpegdupper}.
\end{proof}

\section{Cost of the cut-off}
In this section we estimate the distance between $\rho$ and $\rho_n$. We recall the definition \eqref{deftaun} of $\tnpd$.

\begin{lemma}\label{distgamman}
Let $\rho_n$ be the probability measure defined in \eqref{defgamman}. We have
\begin{equation}\label{statementdistgamman}
nW_p^p(\rho,\rho_n)=o(\tnpd)+ \one_{p>d} O(\tnpd).
\end{equation}
\end{lemma}
\begin{proof}
By subadditivity,
\[
 W_p^p(\rho,\rho_n)\le W_p^p\lt( \one_{(\Omega'_n)^c} \rho, \frac{1-\rho(\Omega'_n)}{\rho(\Omega_n\backslash \Omega'_n)}\one_{\Omega_n\backslash \Omega'_n} \rho\rt).
\]
Let
\[
 \hat{\rho}_n=\frac{1-\rho(\Omega'_n)}{\HH^{d-1}(\partial \Omega'_n)} \HH^{d-1}\restr \partial \Omega'_n.
\]
Then by triangle inequality,
\[
 W_p^p\lt( \one_{(\Omega'_n)^c} \rho, \frac{1-\rho(\Omega'_n)}{\rho(\Omega_n\backslash \Omega'_n)}\one_{\Omega_n\backslash \Omega'_n} \rho\rt)\les W_p^p\lt( \one_{(\Omega'_n)^c} \rho, \hat{\rho}_n\rt)+ W_p^p\lt(\hat{\rho}_n,\frac{1-\rho(\Omega'_n)}{\rho(\Omega_n\backslash \Omega'_n)}\one_{\Omega_n\backslash \Omega'_n} \rho\rt).
\]
We now estimate both terms separately using that since all the measures involved are radially symmetric, so are the optimal transport maps. For the first term we have
\[
 W_p^p\lt( \one_{(\Omega'_n)^c} \rho, \hat{\rho}_n\rt)\les \int_{R'_n}^\infty r^{d-1} |R'_n-r|^p \rho(r)\stackrel{\eqref{condRn}}{\les} \int_{R_n}^\infty r^{d-1} |R_n-r|^p \rho(r).
\]
We then estimate 
\begin{multline*}
\int_{R_n}^{\infty} |R_n-r|^p r^{d-1} \rho(r)dr=  \exp(-V(R_n)) \int_{R_n}^{\infty} |R_n-r|^p r^{d-1} \exp\lt(-\int_{R_n}^r V'(t)dt\rt) dr\\
 \le \rho(R_n) \int_{R_n}^{\infty} |R_n-r|^p r^{d-1} \exp(-C (r^q-R_n^q))
 \les R_n^{d-q+p(1-q)}  \rho(R_n).
\end{multline*}

For the second term we may brutally estimate the transport distance by $R_n-R'_n=\delta_{j_n} R'_n$ (recall \eqref{defrj}) so that
\begin{multline*}
 W_p^p\lt(\hat{\rho}_n,\frac{1-\rho(\Omega'_n)}{\rho(\Omega_n\backslash \Omega'_n)}\one_{\Omega_n\backslash \Omega'_n} \rho\rt)\les (\delta_{j_n} R'_n)^p (1-\rho(\Omega'_n))\stackrel{\eqref{condRn}}{\les} R_n^{p(1-q)}\int_{R_n}^{\infty} r^{d-1}\rho(r)\\
 \les R_n^{d-q+p(1-q)}\rho(R_n).
\end{multline*}

Recalling the definition \eqref{defqalpha} of $\alpha$, $\beta$ as well as \eqref{defbarRn} and \eqref{condRn}, we thus conclude that
\[
 nW_p^p(\rho,\rho_n)\les n R_n^{d-q+p(1-q)}\rho(R_n)\les \begin{cases}
                                                          (\log n)^{\frac{d+p}{q}-1-p} n^{1-\beta} & \textrm{if } p<d\\
                                                          (\log n)^{\frac{d}{q}-p(1-\frac{1}{q})  -1 +\frac{\alpha}{q}} & \textrm{if } p=d\\
                                                          (\log n)^{d-1-p\lt(1-\frac{1}{q}\rt)}  & \textrm{if } p>d.
                                                         \end{cases}
\]
Recalling the definition \eqref{deftaun} of $\tnpd$ and the range \eqref{defqalpha} of $\alpha$ in the case $p=d$ we conclude the proof of \eqref{statementdistgamman}.

\end{proof}

\section{The local terms}
For $r,\delta>0$ we set (recall the definition \eqref{defDdelta} of $D_\delta$)
\[
 \Delta_r^\delta:=r\{ 1\le |x|\le 1+\delta \ : \ \hat{x} \in D_{\delta}\}
\]
where $\hat{x}=x/|x|$. Recall the definitions \eqref{defbetasup} of $\betsup$ and \eqref{defbetinf} of $\betinf$ as well as the definition \eqref{defetanpd} of $\etnpd$. Our main buiding block is a variant of Theorem \ref{upperboundappendix} with a more explicit dependence on the parameters. Let us however stress that the main point in Proposition \ref{prop:localDelta} is that $\Delta^\delta_r$ is almost isometric to a cube while the difficulty in  Theorem \ref{upperboundappendix} is to cover arbitrary domains.
\begin{proposition}\label{prop:localDelta}
 There exists $\omega:\R^+\to \R^+$ with $\lim_{t\to 0}\omega(t)=\lim_{t\to \infty} \omega(t)=0$ with the following property. Let  $r\ge 1$,  $\delta\in (0,1)$  and $\lambda$ be a  Lipschitz continuous probability density on $\Delta_r^\delta$ such that
 \[
  \ell r\delta\le 1,
 \]
 where $\ell:=\|\lambda\|_{\Lip}|\Delta_r^\delta|$.
For every $n \ge 2$, if $(X_i)_{i=1}^n$ are i.i.d. random variables on $\Delta_r^\delta$ with law $\lambda$, denoting  $\mu_n=\sum_{i=1}^n \delta_{X_i}$, we have
 \begin{equation}\label{estimDelta}
  \EE[W^p_{ \Delta_r^\delta}(\mu_n,n \lambda)]
  \le (1+\omega(\ell r\delta)+\omega(\delta)+\omega(n)) \betsup  \etnpd\int_{\Delta_r^\delta} \lambda^{1-\frac{p}{d}}.
 \end{equation}
\end{proposition}
\begin{proof}
 Setting $\lambda_r(x):=r^d\lambda(r x)$  we see by scaling that it is enough to prove \eqref{estimDelta} with $r=1$. To lighten notation we set $\Delta^\delta:=\Delta^\delta_1$. The map
 \begin{eqnarray*}
 \Phi:D_{\delta}\times(0,\delta)\to\mathbb{R}^d,\qquad(z,s)\mapsto (1+s)z\in\Delta^{\delta}
 \end{eqnarray*}
 satisfies $\|\nabla\Phi-{\rm Id}_{\mathbb{R}^d}\|\lesssim\delta$. Moreover, by the very definition of $D_\delta$, the exponential map on $\bS^{d-1}$ provides a diffeomorphism
 $$
 \exp_{e_1}: (-\delta/2, \delta/2)^{d-1} \subseteq\mathbb{R}^{d-1}\to D_{\delta}
 $$
 such that $\|\nabla\exp_{e_1}-{\rm Id}_{\mathbb{R}^d}\|\lesssim\delta$. Letting
 $$
 C_{\delta}:= (-\delta/2, \delta/2)^{d-1} \times(0,\delta),
 $$
 it follows that the map $F:\Delta^{\delta}\mapsto C_\delta$ defined through
 \[
  F^{-1}(y,s):=\Phi(\exp_{e_1}(y),s)=(1+s)\exp_{e_1}(y)
 \]
is a diffeomorphism. Moreover, it satisfies 
\[\|\nabla F-{\rm Id}_{\mathbb{R}^d}\|_\infty +\|\nabla^2 F\|_\infty\lesssim\delta
\]
 and the same estimate  also holds for $F^{-1}$. In particular
 \begin{eqnarray}\label{settori4}
 \left||\Delta^{\delta}|-|C_{\delta}|\right|\leq\int_{C_\delta}|1-\det\nabla F^{-1}|\lesssim\delta|C_{\delta}|.
 \end{eqnarray}
Now we define 
 $$
 \hat{\lambda}:=F_{\sharp}\lambda \qquad \textrm{and } \qquad \hat{\mu}_n:=F_\sharp \mu_n.
 $$
 Notice that $\hat{\mu}_n=\sum_{i=1}^n \delta_{\hat{X}_i}$ where $\hat{X}_i=F(X_i)$ are i.i.d.\ random variables of law $\hat{\lambda}$. We have
 \begin{multline}\label{localhat}
  \hat{\ell}:=\|\hat{\lambda}\|_{\Lip}|C_\delta|\les \ell +\delta \qquad \textrm{and } \\
  \EE[W^p_{ \Delta^\delta}(\mu_n,n \lambda)]\le \|F^{-1}\|_{\Lip}^p \EE[W^p_{ C_\delta}(\hat{\mu}_n,n \hat{\lambda})]\le (1+C\delta)\EE[W^p_{ C_\delta}(\hat{\mu}_n,n \hat{\lambda})].
 \end{multline}
We then set for $x\in C_1$, $\lambda'(x)=\hat{\lambda}(\delta x) \delta^d$ and $\mu_n'=\sum_{i=1}^n \delta_{X_i'}$ where $X_i'=\delta^{-1}\hat{X}_i$ are i.i.d. random variables of law $\lambda'$. By scaling we have 
\begin{multline}\label{localprime}
  \ell':=\|\lambda'\|_{\Lip}|C_1|=\delta \hat{\ell}\les \delta \ell +\delta^2 \quad \textrm{and } \quad
  \EE[W^p_{C_\delta}(\hat{\mu}_n,n \hat{\lambda})]=\delta^p \EE[W^p_{C_1}(\mu'_n,n \lambda')].
 \end{multline}
 Since $C_1$ is a cube, we can use either the Knothe map as in \cite[Lemma 1]{BeCa} or \cite[Proposition 2.4]{ambrosio2022quadratic} to obtain a map $T$ transporting $\lambda'$ on the uniform measure on $C_1$ with
 \[
  \|T\|_{\Lip}, \|T^{-1}\|_{\Lip}\le 1+ O(\ell').
 \]

 Arguing exactly as above and setting $\tilde{\mu}_n=\sum_{i=1}^n \delta_{T(X_i')}$ (so that $T(X_i')$ ae i.i.d. uniformly distributed in $C_1$), we find using the definition of $\etnpd$,
 \begin{multline*}
  \EE[W^p_{C_1}(\mu'_n,n \lambda')]\le (1+O(\ell'))\EE\lt[W^p_{C_1}\lt(\tilde{\mu}_n, \frac{n}{|C_1|}\rt)\rt]
  \le (1+\omega(\ell')+\omega(n)) \betsup  \etnpd |C_1|^{\frac{p}{d}}.
 \end{multline*}

Combining this with \eqref{localhat}, \eqref{localprime} and \eqref{settori4}, we obtain 
\[
  \EE[W^p_{\Delta^\delta}(\mu_n,n \lambda)]
  \le (1+\omega(\delta \ell)+\omega(\delta)+\omega(n)) \betsup  \etnpd|\Delta^\delta|^{\frac{p}{d}}.
  \]
Using that $\lambda |\Delta^\delta|$ is $\ell-$Lipschitz continuous, we can replace $|\Delta^\delta|^{p/d}$ by $\int_{\Delta^\delta} \lambda^{1-p/d}$ at a cost of order $\omega(\delta \ell)$ in the previous inequality thus concluding the proof of \eqref{estimDelta}.
\end{proof}

We recall the following notation from the proof of Proposition \ref{prop:upperbound}. For $\hat{X}_i$ i.i.d. with common law $\rho_n$, we set $\hat{\mu}_n:=\sum_{i=1}^n \delta_{\hat{X}_i}$ and then
\[\kappa_{1}:= \frac{\hat{\mu}_n(B_{r_1})}{\rho_n(B_{r_1})} \qquad \textrm{and } \qquad \kappa_j^z:= \frac{\hat{\mu}_n(\Delta_j^z)}{\rho_n(\Delta_j^z)}.\] 
\begin{proposition}\label{prop:localterms}
 We have for  $p\le d$
 \begin{multline}\label{eq:proplocaltermssmall}
  \EE[W_{B_{r_1}}^p(\hat{\mu}_n, \kappa_1 \rho_n)] +\sum_{j=2}^{j_n} \frac{1}{\stilde_j } \int_{\bStilde^{d-1}} \EE[W_{\Delta_j^z}^p(\hat{\mu}_n, \kappa_j^z \rho_n)] d\stilde(z)\\
  \le \lt(1+\omega(\eps)+\omega(n)\rt) \betsup \etnpd\int_{\Omega_n} \rho_n^{1-\frac{p}{d}} \lt(1+\one_{d=2} \lt(\frac{\log n \rho_n}{\log n}\rt)^{\frac{p}{2}}\rt)
 \end{multline}
and for $p>d$
\begin{equation}\label{eq:proplocaltermslarge}
  \EE[W_{B_{r_1}}^p(\hat{\mu}_n, \kappa_1 \rho_n)] +\sum_{j=2}^{j_n} \frac{1}{\stilde_j} \int_{\bStilde^{d-1}} \EE[W_{\Delta_j^z}^p(\hat{\mu}_n, \kappa_j^z \rho_n)] d\stilde(z)\les_\eps n^{1-\frac{p}{d}}\int_{\Omega_n} \rho_n^{1-\frac{p}{d}}.
 \end{equation}
\end{proposition}
\begin{proof}
  We set
 \[n_1:=\hat{\mu}_n(B_{r_1}) \qquad \textrm{ and }  \qquad n_j^z:= \hat{\mu}_n(\Delta_j^z).\]
In $B_{r_1}$ we use Theorem \ref{th:upperboundappendix} to obtain
\begin{multline*}
 \EE[W_{B_{r_1}}^p(\hat{\mu}_n, \kappa_1 \rho_n)]=\EE[W_{B_{r_1}}^p(\hat{\mu}_n, n_1 \frac{\rho}{\rho(B_{r_1})})]\\
 \le \EE\lt[(1+\omega(n_1))\lt(\frac{n_1}{\rho(B_{r_1})}\rt)^{1-\frac{p}{d}} \lt(1+\one_{d=2} (\log n_1)^{\frac{p}{2}}\rt)\rt]\betsup\int_{B_{r_1}} \rho_n^{1-\frac{p}{d}}.
\end{multline*}
Since $n_1$ is a binomial random variable with exponential concentration around its expectation $\EE[n_1]=n\rho(B_{r_1})$,
\begin{multline}\label{localball}
 \EE[W_{B_{r_1}}^p(\hat{\mu}_n, \kappa_1 \rho_n)]\le (1+\omega(n))\betsup\etnpd\int_{B_{r_1}} \rho_n^{1-\frac{p}{d}}\\
 \le \lt(1+\omega(n)\rt) \betsup n^{1-\frac{p}{d}}\int_{B_{r_1}} \rho_n^{1-\frac{p}{d}} \lt(1+\one_{d=2} (\log n \rho_n)^{\frac{p}{2}}\rt).
\end{multline}
We now consider a fixed set $\Delta^z_j$. Since $\lambda_j^z:=\rho_n/\rho_n(\Delta_j^z)$ is Lipschitz continous in $\Delta^z_j$ with
\[
 \ell_j^z:=\|\lambda_j^z\|_{\Lip}|\Delta_j^z|\les \sup_{\Delta_j^z} |\nabla V|\les r_j^{q-1}
\]
we find 
\begin{equation}\label{lipgammaDeljz}
 \ell_j^z r_j \delta_j\les r_j^q \delta_j\stackrel{\eqref{defdelta}\&\eqref{condRn}}\les \eps.
\end{equation}
We may thus apply \eqref{estimDelta} from Proposition \ref{prop:localDelta} and obtain
\begin{multline*}
 \EE[W_{\Delta_j^z}^p(\hat{\mu}_n, \kappa_j^z \rho_n)]\le\\
  \betsup  \EE\lt[(1+\omega(\eps)+\omega(n_j^z))\lt(\frac{n_j^z}{\rho_n(\Delta_j^z)}\rt)^{1-\frac{p}{d}}(1+\one_{d=2}(\log (1+n_j^z))^{\frac{p}{2}}) \Bigg| n_j^z\ge 1 \rt]\int_{\Delta_j^z} \rho_n^{1-\frac{p}{d}}.
\end{multline*}

We now observe that $n_j^z$ is a binomial random variable with 
\begin{equation}\label{Enjz}
 \EE[n_j^z]=n \rho_n(\Delta_j^z)\simeq n\rho_n(r_j) (r_j\delta_j)^d\stackrel{\eqref{defdelta}}{\ges_\eps}  n\rho(R_n) R_n^{d(1-q)}.
\end{equation}
From \eqref{condRn} and \eqref{defbarRn} we get (recall the definition \eqref{defqalpha} of $\alpha$ and $\beta$)
\[
 n\rho(R_n) R_n^{d(1-q)}\simeq \begin{cases}
                            n^{1-\beta} (\log n)^{\frac{d}{q}(1-q)} &\textrm{if } p<d\\
                            (\log n)^{\frac{1}{q}(\alpha-d(q-1))} &\textrm{if } p=d\\
                            1 &\textrm{if } p>d.
                           \end{cases}
\]
We now argue separately in the cases $p\le d$ and $p>d$. The case $d=2$ requires a little more care so we first consider the case $d\ge 3$. Since $p\le d$, we have   $\lim_{n\to \infty} \inf_{j,z} \EE[n_j^z]=\infty$ and by the good concentration properties of binomial random variables we have 
\[
 \EE\lt[(1+\omega(\eps)+\omega(n_j^z))\lt(\frac{n_j^z}{\rho(\Delta_j^z)}\rt)^{1-\frac{p}{d}} \Bigg| n_j^z\ge 1\rt]
 \le (1+\omega(\eps)+\omega_\eps(n))n^{1-\frac{p}{d}}.
\]
This proves that for $p\le d$,
\begin{equation}\label{estimlocp<d}
 \EE[W_{\Delta_j^z}^p(\hat{\mu}_n, \kappa_j^z \rho_n)]\le \betsup (1+\omega(\eps)+\omega_\eps(n))n^{1-\frac{p}{d}}\int_{\Delta_j^z} \rho_n^{1-\frac{p}{d}}.
\end{equation}
After summation and integration, combining \eqref{localball} and \eqref{estimlocp<d} concludes the proof of \eqref{eq:proplocaltermssmall}. When $p<d=2$,  since $n \rho_n(\Delta_j^z)= C_\eps n^{1-\beta}$ for $\Delta_j^z$ close to $\partial B_{R_n}$ we need to be a bit more careful. Letting $\hat{R}_n:=\log\log n$ (this choice is somewhat arbitrary) we have for $\Delta_j^z\subset  B_{\hat{R}_n}$,  $\log (n \rho_n(\Delta_j^z))= (1+\omega_\eps(n))\log n$ while for $\Delta_j^z\cap A_{\hat{R}_n,R_n}\neq \emptyset$ (recall definition \eqref{defannulus}), $\log (n \rho_n(\Delta_j^z))=O_\eps(\log n)$ so that the same argument as above gives
\begin{multline*}
 \EE[W_{B_{r_1}}^p(\hat{\mu}_n, \kappa_1 \rho_n)] +\sum_{j=2}^{j_n} \frac{1}{\stilde_j} \int_{\bStilde^{1}} \EE[W_{\Delta_j^z}^p(\hat{\mu}_n, \kappa_j^z \rho_n)] d\stilde(z)\\
  \le \lt(1+\omega(\eps)+\omega(n)\rt) \betsup n^{1-\frac{p}{2}} (\log n)^{\frac{p}{2}}\lt(\int_{B_{\hat{R}_n}} \rho_n^{1-\frac{p}{2}} +C_\eps \int_{A_{\hat{R}_n,R_n}}\rho_n^{1-\frac{p}{2}}\rt)\\
  \le \lt(1+\omega(\eps)+\omega(n)\rt) \betsup n^{1-\frac{p}{2}} \int_{\Omega_n} \rho_n^{1-\frac{p}{2}}(\log n \rho_n)^{\frac{p}{2}}.
\end{multline*}
Thus \eqref{eq:proplocaltermssmall} also holds in this case. If $p=d=2$ we have \[
    \EE[W_{\Delta_j^z}^2(\hat{\mu}_n, \kappa_j^z \rho_n)]\le \betsup (1+\omega(\eps)+\omega_\eps(n))\int_{\Delta_j^z} \log(n \rho_n |\Delta_j^z|).
                                                                            \]
Since it is readily checked that
\[
 \sum_{j=2}^{j_n} \frac{1}{\stilde_j} \int_{\bStilde^{1}} |\Delta_j^z| \log |\Delta_j^z| d\stilde(z)=o(\etnpd),
\]
after summation we conclude again the proof of \eqref{eq:proplocaltermssmall}.
\\

In the case $p>d$, since $\inf_{j,z} \EE[n_j^z]\ges_\eps 1$, we can still conclude that 
\begin{multline*}
 \EE\lt[(1+\omega(\eps)+\omega(n_j^z))\lt(\frac{n_j^z}{\rho(\Delta_j^z)}\rt)^{1-\frac{p}{d}}(1+\one_{d=2}(\log (1+n_j^z))^{\frac{p}{2}}) \Bigg| n_j^z\ge 1 \rt]\\
 \les_\eps n^{1-\frac{p}{d}}(1+\one_{d=2}(\log( 1+n \rho_n(\Delta_j^z)))^{\frac{p}{2}}).
\end{multline*}
Therefore 
\begin{equation}\label{estimlocp>d}
 \EE[W_{\Delta_j^z}^p(\hat{\mu}_n, \kappa_j^z \rho_n)]\les_\eps n^{1-\frac{p}{d}}(1+\one_{d=2}(\log( 1+n \rho_n(\Delta_j^z)))^{\frac{p}{2}})\int_{\Delta_j^z} \rho_n^{1-\frac{p}{d}}.
 \end{equation}
In the case $d>2$, after summation and integration, combining \eqref{localball} and \eqref{estimlocp>d} concludes the proof of \eqref{eq:proplocaltermslarge}. We are left with the case $p>d=2$. Using that $\rho_n(\Delta_j^z)\simeq_\eps  \rho(r_j) r_j^{d} (\log n)^{-d}$ we have
\begin{multline*}
 \sum_{j=2}^{j_n} \frac{1}{\stilde_j} \int_{\bStilde^{d-1}}(\log( 1+n \rho_n(\Delta_j^z)))^{\frac{p}{2}}\int_{\Delta_j^z} \rho_n^{1-\frac{p}{d}} d\stilde(z) \\\les_\eps \int_{\Omega_n\backslash B_{r_1}} \lt(\log\lt( 1+n \lt(\frac{|x|}{\log n}\rt)^{d}\rho_n \rt)\rt)^{\frac{p}{2}} \rho_n^{1-\frac{p}{d}}.
\end{multline*}
Before starting the estimate let us recall that by \eqref{estim:integrhonp>d}
\begin{equation}\label{integp>d}
 \int_{\Omega_n} \rho_n^{1-\frac{p}{d}}\simeq_\eps  \rho(R_n)^{1-\frac{p}{d}}R_n^{d-q}.
\end{equation}
We first estimate in $B_{R_n/2}\backslash B_{r_1}$,
\begin{multline*}
 \int_{B_{R_n/2}\backslash B_{r_1}} \lt(\log\lt( 1+n \lt(\frac{|x|}{\log n}\rt)^{d}\rho_n \rt)\rt)^{\frac{p}{2}} \rho_n^{1-\frac{p}{d}}\les (\log n)^{\frac{p}{2}} \int_{B_{R_n/2}} \rho^{1-\frac{p}{d}}\\
 \les (\log n)^{\frac{p}{2}} R_n^{d-q} \rho(R_n/2)^{1-\frac{p}{d}}\les \rho(R_n)^{1-\frac{p}{d}}R_n^{d-q} \lt((\log n)^{\frac{p}{2}} \exp(-C R_n^q)\rt) \\\stackrel{\eqref{integp>d}}{\ll}\int_{\Omega_n} \rho_n^{1-\frac{p}{d}}.
\end{multline*}
Here we used that by \eqref{hypV}, we have for some $C>0$,
\[
 V(R_n/2)\le V(R_n)-C R_n^q.
\]

In $\Omega_n\backslash B_{R_n/2}$ we use that since $n \lt(|x|/\log n\rt)^{d}\rho_n \ges_\eps 1$ we have for some $0<c_\eps\ll 1$,
\[\lt( 1+n \lt(\frac{|x|}{\log n}\rt)^{d}\rho_n \rt)\le 1+ c_\eps \log \lt(n \lt(\frac{|x|}{\log n}\rt)^{d}\rho_n\rt)\les 1+c_\eps\log (n R_n^{d(1-q)}\rho_n).\]
Since $ nR_n^{d(1-q)} \rho(R_n)\simeq_\eps 1$ we can write
\[
 n R_n^{d(1-q)}\rho_n(x)\simeq_\eps \frac{\rho_n(x)}{\rho_n(R_n)}=\exp(V(R_n)-V(x)).
\]
Thus
\begin{multline*}
 \int_{\Omega_n\backslash B_{R_n/2}} \lt(\log\lt( 1+n \lt(\frac{|x|}{\log n}\rt)^{d}\rho_n \rt)\rt)^{\frac{p}{2}} \rho_n^{1-\frac{p}{d}}\les_\eps \int_{\Omega_n\backslash B_{R_n/2}} (1+(V(R_n)-V(x)))^{\frac{p}{2}}\rho_n^{1-\frac{p}{d}}\\
 \stackrel{\eqref{hypV}}{\les} R_n^{d-1} \int_{R_n/2}^{R_n} (1+  C  (R_n^q-r^q))^{\frac{p}{2}} \rho_n^{1-\frac{p}{d}}(r)
 \les R_n^{d-q} \rho_n^{1-\frac{p}{d}}(R_n)\\
 \stackrel{\eqref{integp>d}}{\les}\int_{\Omega_n} \rho_n^{1-\frac{p}{d}}.
\end{multline*}
In conclusion we find
\[
 \int_{\Omega_n\backslash B_{r_1}} \lt(\log\lt( 1+n \lt(\frac{|x|}{\log n}\rt)^{d}\rho_n \rt)\rt)^{\frac{p}{2}} \rho_n^{1-\frac{p}{d}}\les_\eps \int_{\Omega_n} \rho_n^{1-\frac{p}{d}}.
\]
Since by \eqref{integp>d} again,
\[
 (\log n)^{\frac{p}{2}}\int_{B_{r_1}} \rho_n^{1-\frac{p}{d}}\ll \int_{\Omega_n} \rho_n^{1-\frac{p}{d}},
\]
this concludes the proof of \eqref{eq:proplocaltermslarge} also in the case $d=2$.
\end{proof}

\section{The global terms}
\subsection{The main radial estimates}\label{sec:radial}
\begin{definition}
 Let $\gamma\in \R$. We say that a probability distribution $\lambda$ on $(0,R)$ is of type $\gamma$ if there exist $\bar R\in[0,R]$ such that  there exist  $C_{\pm}$ positive with
 \begin{equation}\label{cmoins}
  \int_0^r  \lambda \le C_- r \lambda(r) \qquad \textrm{for } r\in (0,\bar R)
 \end{equation}
and
\begin{equation}\label{cplus}
 \int_r^R \lambda \le C_+ r^{\gamma}\lambda(r) \qquad \textrm{for } r\in (\bar R, R).
\end{equation}
We write $A\les_\pm B$ if there is a constant $C>0$ depending only on $p,d,\gamma$ and $C_{\pm}$ such that $A\le C B$.
\end{definition}
We now prove two very similar result. It should be possible to combine them in a single statement but we decided for clarity to keep them separate.
\begin{proposition}\label{prop:radial1}
 Let $\lambda$ be of  type $\gamma$ and assume that $(0,R)=(0,\bar R)\cup [\bar R, R)=I_0\cup I_1$. For $n\in \N$, let $(n_i)_{i\in \{0,1\}}$ be a multinomial distribution of parameters $(n \lambda(I_i))_{i\in \{0,1\}}$ and set $\kappa_i=n_i/\lambda(I_i)$. Assume that
 \begin{equation}\label{hyp:multinomial}
\quad n\lambda(I_i)\ge 1 \qquad \textrm{for } i=0,1.
 \end{equation}
 Then,
\begin{equation}\label{estimradialabstract}
 \EE[W_p^p(\sum_i \kappa_i  \lambda \one_{I_i},n \lambda)]\les_\pm n^{1-\frac{p}{2}} \lt(\bar{R}^p \lambda(I_0)^{1-\frac{p}{2}}+ \lambda(I_1)^{-\frac{p}{2}} \int_{I_1} r^{\gamma p} \lambda(r)\rt).
\end{equation}

\end{proposition}
\begin{proof}

Let $\phi$ be the solution to
\[
 \phi''=\sum_i (\kappa_i-n)  \lambda \one_{I_i}, \qquad \phi'(0)=\phi'(R)=0
\]
we then have by Proposition \ref{BBcomodo},
\[
 \EE[W_p^p(\sum_i \kappa_i  \lambda \one_{I_i}, n\lambda)]\les  n^{1-p}\int_{0}^R \frac{1}{\lambda^{p-1}}\EE[|\phi'|^p].
\]
We claim that
\begin{equation}\label{claimcmoins}
  n^{1-p}\int_{I_0} \frac{1}{\lambda^{p-1}}\EE[|\phi'|^p]\les_\pm n^{1-\frac{p}{2}} \bar{R}^p \lambda(I_0)^{1-\frac{p}{2}}.
\end{equation}
  We write for $r\in I_0$,
\[
 |\phi'(r)|=\lt|\frac{n_0-n\lambda(I_0)}{\lambda(I_0)} \int_0^r \lambda\rt|\stackrel{\eqref{cmoins}}{\les_\pm} \frac{|n_0-n\lambda(I_0)|}{\lambda(I_0)} r \lambda(r)
\]
so that using \eqref{hyp:multinomial} we get
\[
 \EE[|\phi'(r)|^p]\les_\pm n^{\frac{p}{2}} \lambda(I_0)^{-\frac{p}{2}} r^p \lambda(r)^p.
\]
After integration this yields
\[
 n^{1-p}\int_0^{\bar R} \frac{1}{\lambda^{p-1}}\EE[|\phi'|^p]\les_\pm n^{1-\frac{p}{2}} \lambda(I_0)^{-\frac{p}{2}}\int_0^{\bar R} r^p \lambda\les_\pm n^{1-\frac{p}{2}} \bar{R}^p \lambda(I_0)^{1-\frac{p}{2}}.
\]
This proves \eqref{claimcmoins}. We now claim that
\begin{equation}\label{claimcplus}
n^{1-p}\int_{I_1} \frac{1}{\lambda^{p-1}}\EE[|\phi'|^p]\les_\pm  n^{1-\frac{p}{2}}  \lambda(I_1)^{-\frac{p}{2}} \int_{I_1} r^{\gamma p} \lambda(r).
\end{equation}
Let $r\in I_1=[\bar R,R)$. We write
\[
 |\phi'(r)|=\lt|\int_r^{R} (\kappa_1-n)\lambda\rt|= \frac{|n_1 -n\lambda(I_1)|}{\lambda(I_1)}\int_{r}^{ R} \lambda.
\]
Therefore, using \eqref{hyp:multinomial} we get
\[
 \EE[|\phi'|^p]\les_{\pm}  n^{\frac{p}{2}}\lambda(I_1)^{-\frac{p}{2}} \lt(\int_{r}^{ R} \lambda\rt)^p\stackrel{\eqref{cplus}}{\les}_{\pm} n^{\frac{p}{2}}\lambda(I_1)^{-\frac{p}{2}} r^{\gamma p} \lambda(r)^p.
\]
After integration this proves \eqref{claimcplus}. Combining \eqref{claimcmoins} and \eqref{claimcplus} together proves \eqref{estimradialabstract}.
\end{proof}
We now prove a second version of this estimate when $(\bar R,R)$ is decomposed into sub-intervals.
\begin{proposition}\label{prop:radial2}
 Let $\lambda$ be of  type $\gamma$  with $\lambda=0$ on $(0,\bar R)$. Assume that $(\bar R,R)=\cup I_i$ where $I_i$ are disjoint intervals indexed in increasing order. For $n\in \N$, let $(n_i)_{i}$ be a multinomial distribution of parameters $(n \lambda(I_i))_{i}$ and set $\kappa_i=n_i/\lambda(I_i)$. Assume that for some $\eps>0$,
 \begin{equation}\label{hyp:multinomial2}
\quad n\lambda(I_i)\ge \eps \qquad \textrm{for  every } i\ge1.
 \end{equation}
 Assume moreover that $r\mapsto r^{\gamma} \lambda(r)$ is decreasing in each $I_i$. Then,
\begin{equation}\label{estimradialabstract2}
 \EE[W_p^p(\sum_i \kappa_i  \lambda \one_{I_i},n \lambda)]\les_{\pm,\eps} n^{1-\frac{p}{2}}\int_{\bar R}^R r^{\gamma \frac{p}{2}} \lambda(r)^{1-\frac{p}{2}} dr.
\end{equation}

\end{proposition}
\begin{proof}
 As above we consider  $\phi$  the solution to
\[
 \phi''=\sum_{i} (\kappa_i-n)  \lambda \one_{I_i}, \qquad \phi'(\bar R)=\phi'(R)=0
\]
so that as above,
\[
 \EE[W_p^p(\sum_i \kappa_i  \lambda \one_{I_i}, n\lambda)]\les  n^{1-p}\int_{\bar R}^R \frac{1}{\lambda^{p-1}}\EE[|\phi'|^p].
\]
Let $r\in I_i=(r_{i-1},r_i)$ for some $i\ge 1$. We write
\begin{equation}\label{eq:phi'plus}
 \phi'(r)=-\int_r^{r_i} (\kappa_i-n)\lambda -\int_{r_i}^R \sum_{j>i} (\kappa_j-n)\lambda.
\end{equation}
Regarding the first term we have
\[
 \lt|\int_r^{r_i} (\kappa_i-n)\lambda\rt|= \frac{|n_i -n\lambda(I_i)|}{\lambda(I_i)}\int_{r}^{r_i} \lambda.\]
For the second term we observe that
\[
 \int_{r_i}^R \sum_{j>i} (\kappa_j-n)\lambda=\sum_{j>i} \lambda(I_j)(\kappa_j-n)=\lt( \sum_{j>i} n_j\rt)- n \lt(\sum_{j>i} \lambda(I_j)\rt)= n_i^+-n \lambda(I_i^+)
\]
where  $n_i^+=\sum_{j>i} n_j$ and $I_i^+=(r_i,R)=\cup_{j>i} I_j$ so that $\sum_{j>i} \lambda(I_j)=\lambda(I_i^+)$.  Observe that  $n_i^+$ is a binomial random variable of parameters $(n, \lambda(I_i^+))$. Plugging this into \eqref{eq:phi'plus} we find thanks to \eqref{hyp:multinomial2},
\[
 \EE[|\phi'|^p]\les \EE[|n_i-n\lambda(I_i)|^p]+\EE[|n_i^+-n \lambda(I_i^+)|^p]\les_\eps n^{\frac{p}{2}}\lambda(I_i)^{-\frac{p}{2}} \lt(\int_{r}^{r_i} \lambda\rt)^p + n^{\frac{p}{2}} \lambda(I_i^+)^{\frac{p}{2}}.
 \]
 Using that
 \[
  \lambda(I_i)^{-\frac{p}{2}} \lt(\int_{r}^{r_i} \lambda\rt)^p\le \lt(\int_{r}^{r_i} \lambda\rt)^{\frac{p}{2}}\le \lambda(I_i^+)^{\frac{p}{2}}
 \]
this simplifies to
 \[
 \EE[|\phi'|^p]\les_\eps n^{\frac{p}{2}} \lambda(I_i^+)^{\frac{p}{2}} \stackrel{\eqref{cplus}}{\les_{\pm}} n^{\frac{p}{2}} r_i^{\gamma\frac{p}{2}} \lambda(r_i)^{\frac{p}{2}}\le  n^{\frac{p}{2}} r^{\gamma\frac{p}{2}} \lambda(r)^{\frac{p}{2}}
\]
where in the last inequality we used that $r\mapsto r^{\gamma} \lambda(r)$ is decreasing in $I_i$. After integration this gives
\[
 n^{1-p}\int_{\bar R}^R \frac{1}{\lambda^{p-1}}\EE[|\phi'|^p]\les_{\pm,\eps} n^{1-\frac{p}{2}}\int_{\bar R}^R r^{\gamma\frac{p}{2}} \lambda(r)^{1-\frac{p}{2}}
\]
concluding the proof of \eqref{estimradialabstract2}.
\end{proof}
\subsection{The main angular estimate}\label{sec:angular}
Recalling the definition \eqref{defLambda} of $\Lambda_{r,R}^{z,\delta}$ we set $\Lambda^{z,\delta}:=\Lambda_{0,\infty}^{z,\delta}$.
\begin{proposition}\label{prop:angle}
 Let $\lambda$ be a radially symmetric probability density on $\R^d$, $n\in \N$ and $\delta\in(0,\pi]$. Let
 \[
  M_p:=\int_{\R^d} |x|^p d\lambda \quad \text{and} \quad \stilde_\delta := \int_{\bStilde^{d-1}} \one_{D_{\delta}(z)}(e_1) d \stilde(z).
 \]
For $(X_i)_{i=1}^n$ i.i.d. random variables of law $\lambda$, let $\mu_n:=\sum_{i=1}^n \delta_{X_i}$ and for $z \in \bStilde^{d-1}$
 \[
 \kappa_z:=\frac{\mu_n(\Lambda^{z,\delta})}{\lambda(\Lambda^{z,\delta})}.
 \]
Then, it holds (recall the definition \eqref{defetanpd} of $\etnpd$)
\begin{equation}\label{eq:angle}
  \EE\lt[ W_p^p\lt( \frac{1}{\stilde_\delta }\int_{\bStilde^{d-1}} \kappa_z \one_{\Lambda^{z,\delta}} \lambda d\stilde(z), n\lambda \rt)\rt]\\
  \les M_p \,\eta_n^{p,d-1}.
\end{equation}
\end{proposition}

\begin{proof}
We use polar coordinates to write for $r>0$ and $z = (y,U)\in \bStilde^{d-1}$,
\[\lambda(r y)=\lambda_{\rm rad}(r) dr \otimes d\sigma(y).\]
Noticing that for $y' \in \bS^{d-1}$ it holds $\one_{\Lambda^{z,\delta}}(r y')=\one_{D_\delta(z)}(y')$, we can write
\[
 \lt(\frac{1}{\stilde_d}\int_{\bStilde^{d-1}} \kappa_z \one_{\Lambda^{z,\delta}} \lambda d\stilde(z)\rt)(r y')=\lambda_{\rm rad}(r) dr \otimes \nu(y')d\sigma(y')
\]
where
\[
 \nu(y'):=\frac{1}{\stilde_\delta}\int_{\bStilde^{d-1}} \kappa_z \one_{D_\delta(z)}(y') d\stilde(z).
\]
If $T_{\rm ang}$ is the optimal transport map between $n\sigma$ and $\nu\sigma$ (for the Euclidean cost in $\R^d$), we may use  $T(ry):= rT_{\rm ang}(y)$ as a competitor to estimate
\begin{multline}\label{estimang1}
 W_p^p\lt( \frac{1}{\stilde_\delta}\int_{\bStilde^{d-1}} \kappa_z \one_{\Lambda^{z,\delta}} \lambda d\stilde(z), n\lambda \rt)\le \int_0^\infty \int_{\bS^{d-1}} |T(r y)- ry|^p \lambda_{\rm rad}(r) dr n d\sigma(y)\\
 =\int_0^\infty r^p  \lambda_{\rm rad}(r) dr\int_{\bS^{d-1}} |T_{\rm ang}( y)- y|^p nd\sigma(y)\\
 =M_p W_p^p\lt(\nu\sigma, n\sigma\rt)\le M_p W_{\bS^{d-1}}^p\lt(\nu\sigma, n\sigma\rt)
\end{multline}
where with a slight abuse of notation we write $W_{\bStilde^{d-1}}^p$ for the Wasserstein distance on the sphere with respect to the geodesic distance $\Sdist$.
Let now $\hat{X}_i:=X_i/|X_i|$ be the projections of $X_i$ on the sphere. Since $\lambda$ is radially symmetric, $\hat{X}_i$ are i.i.d.\ uniformly distributed on $\bS^{d-1}$ i.e. their common law is $\sigma$. If $\hat{\mu}_n:=\sum_{i=1}^n \delta_{\hat{X}_i}$, using that $\lambda(\Lambda^{z,\delta})=\stilde_j$, we have
\begin{equation}\label{eq:kappazhatmun}
 \kappa_z=\frac{\hat{\mu}_n(D_\delta(z))}{\sigma(D_\delta)}.
\end{equation}
We define the  Markov kernel $k_{\delta}$ on $\mathbb{S}^{d-1}\times\mathbb{S}^{d-1}$ as
\begin{eqnarray}\label{globalesettori5}
k_{\delta}(x,y):= \int_{\bStilde^{d-1}} \frac{ \mathbbm{1}_{D_{\delta}(z)}(x)\mathbbm{1}_{D_{\delta}(z)}(y)}{\stilde_\delta \sigma(D_\delta)} d\stilde(z).
\end{eqnarray}
Given a measure $\mu$ on $\mathbb{S}^{d-1}$, we then define $K_{\delta}(\mu)$ as the measure on $\mathbb{S}^{d-1}$ whose density with respect to $\sigma$ is $\int_{\mathbb{S}^{d-1}}k_{\delta}(x,\cdot)d\mu$, i.e.
\begin{eqnarray}\label{globalesettori6}
K_{\delta}(\mu)(E):=\int_{\mathbb{S}^{d-1}}\mathbbm{1}_E(x)\int_{\mathbb{S}^{d-1}}k_{\delta}(x,y)d\mu(y)d\sigma(x).
\end{eqnarray}
In particular, by \eqref{eq:kappazhatmun} we have
\begin{eqnarray}\label{globalesettori7}
\nu\sigma=K_{\delta}\left(\hat{\mu}_n\right)\qquad \textrm{and } \qquad n\sigma=K_{\delta}\left( n\sigma\right).
\end{eqnarray}
For $x,y\in\bS^{d-1}$, let   $\eta_y(x):=k_\delta(x,y)$. We claim (and will prove below) that for every $y,y'\in \bS^{d-1}$  we have (identifying $\eta_y \sigma$ with $\eta_y$)
\begin{equation}\label{claimWKdelta}
 W_{\bS^{d-1}}^p(\eta_y,\eta_{y'})\le \Sdist^p(y,y').
\end{equation}
Postponing the proof of \eqref{claimWKdelta} we now estimate $W_{\bS^{d-1}}^p(K_{\delta}(\hat{\mu}_n),K_{\delta}( n\sigma))$. Let $\hat{\pi}$ be an optimal transport plan for $W_{\bS^{d-1}}^p(\hat{\mu}_n, n\sigma)$ and for $y,y'\in \bS^{d-1}$, let $\pi_{y,y'}$ be an optimal transport plan for $W_{\bS^{d-1}}^p(\eta_y,\eta_{y'})$ we define the transport plan $\pi$ by its action on $\zeta\in C^0(\bS^{d-1}\times \bS^{d-1})$:
\[
 \int_{\bS^{d-1}\times\bS^{d-1}} \zeta(x,x')d\pi(x,x'):=\int_{\bS^{d-1}\times\bS^{d-1}}\lt[\int_{\bS^{d-1}\times\bS^{d-1}} \zeta(x,x') d\pi_{y,y'}(x,x')\rt]d\hat{\pi}(y,y').
\]
Let us prove that the first marginal of $\pi$ is $K_{\delta}(\hat{\mu}_n)$. Indeed, for every $\zeta\in C^0(\bS^{d-1})$, we have
\begin{multline*}
 \int_{\bS^{d-1}\times\bS^{d-1}} \zeta(x)d\pi(x,x')=\int_{\bS^{d-1}\times\bS^{d-1}}\lt[\int_{\bS^{d-1}\times\bS^{d-1}} \zeta(x) d\pi_{y,y'}(x,x')\rt]d\hat{\pi}(y,y')\\
 =\int_{\bS^{d-1}\times\bS^{d-1}}\lt[\int_{\bS^{d-1}} \zeta(x) \eta_y(x) d\sigma(x)\rt]d\hat{\pi}(y,y')
 =\int_{\bS^{d-1}}\lt[\int_{\bS^{d-1}} \zeta(x) k_\delta(x,y) d\sigma(x)\rt]d\hat{\mu}_n(y)\\
 =\int_{\bS^{d-1}} \zeta(x) d K_\delta(\hat{\mu}_n)(x).
\end{multline*}
Similarly the second marginal of $\pi$ is $K_{\delta}( n\sigma)$. Therefore,
\begin{align*}
 W_{\bS^{d-1}}^p&\lt(K_\delta\lt(\hat{\mu}_n\rt), K_\delta\lt(n\sigma\rt)\rt)
 \\
 &\le
\int_{\bS^{d-1}\times\bS^{d-1}}\lt[\int_{\bS^{d-1}\times\bS^{d-1}} \Sdist^p(x,x') d\pi_{y,y'}(x,x')\rt]d\hat{\pi}(y,y')\\
&=\int_{\bS^{d-1}\times\bS^{d-1}}W_{\bS^{d-1}}^p(\eta_y,\eta_{y'})d\hat{\pi}(y,y') \\&\stackrel{\eqref{claimWKdelta}}{\le}\int_{\bS^{d-1}\times\bS^{d-1}}\Sdist^p(y,y')d\hat{\pi}(y,y')\\
&=W_{\bS^{d-1}}^p\lt(\hat{\mu}_n, n\sigma\rt).
 \end{align*}
Combining this with \eqref{estimang1},  \eqref{globalesettori7} together with the scaling of the optimal matching problem on the sphere, see e.g. \cite[Section 4]{Le17}, we conclude that
\begin{multline*}
 \EE\lt[W_p^p\lt( \frac{1}{\stilde_j}\int_{\bStilde^{d-1}} \kappa_z \one_{\Lambda^{z,\delta}} \lambda d\stilde(z), n\lambda \rt)\rt]\les M_p \EE\lt[W_{\bS^{d-1}}^p\lt(\hat{\mu}_n, n\sigma\rt)\rt]
 \les M_p\,\eta_n^{p,d-1}.
\end{multline*}
Therefore in order to conclude the proof of \eqref{eq:angle} we are left with the proof of \eqref{claimWKdelta}. Let $y,y'\in \bS^{d-1}$ and let $\Pi$ be the two-dimensional plane containing the great circle of $\bS^{d-1}$ containing $y$ and $y'$. We then let $R_\theta$ be the rotation of angle $\theta\in(0,\pi]$ such that $R_\theta={\rm Id}$ on $\Pi^\perp$ and $R_\theta(y)=y'$ so that $\Sdist(y,y')=\theta$. Since $R_\theta$ is an isometry of $\bS^{d-1}$, it induces an orthogonal transformation $U_{\theta}$ between the tangents at $y$ and $y'$, i.e., we can write
\[ R_\theta \left( D_{\delta}(y, U) \right)  = D_{\delta}(y', U_\theta U)\]
Since in the definition of $\eta_y(x) = k_\delta(x,y)$ we integrate (also) with respect to the Haar measure on $O(d-1)$, which is left invariant, we have that $(R_\theta)_\sharp \eta_y = \eta_{y'}$. Finally, for every $x\in \bS^{d-1}$ if $r\le 1$ is the radius of the circle given by the intersection between $\bS^{d-1}$ and the plane parallel to $\Pi$ containing $x$ we have
\[
 \Sdist(x,R_\theta(x))\le r \theta\le \Sdist(y,y')
\]
so that
\[
 W_{\bStilde^{d-1}}^p(\eta_y,\eta_{y'})\le \int_{\bStilde^{d-1}} \Sdist^p(x,R_\theta(x)) \eta_y(x) d\sigma(x)\le \Sdist^p(y,y').
\]
This concludes the proof of  \eqref{claimWKdelta}.
\end{proof}

\subsection{Estimate of the global terms}
We may now estimate the global cost. Let us recall that if $\hat{\mu}_n=\sum_{i=1}^n\delta_{\hat{X}_i}$ with $\hat{X}_i$ i.i.d. random variables with law $\rho_n$, we have set in \eqref{defkappajz} (recall the definition \eqref{defDeltazj} of $\Dzj$),
\[\kappa_{1}:= \frac{\hat{\mu}_n(B_{r_1})}{\rho_n(B_{r_1})} \qquad \textrm{and } \qquad \kappa_j^z:= \frac{\hat{\mu}_n(\Dzj)}{\rho_n(\Dzj)}.\]

We recall also the definition \eqref{deftaun} of $\tnpd$.
\begin{proposition}\label{prop:global}
We have
\begin{equation}\label{eq:global}
 \EE\lt[W_p^p\lt(\kappa_1 \one_{B_{r_1}}\rho_n+\sum_{j=2}^{j_n}\frac{1}{\stilde_j }\int_{\bStilde^{d-1}}\kappa_j^z \one_{\Delta_j^z} \rho_n d\stilde(z), n\rho_n\rt)\rt]=o_\eps(\tnpd)+\one_{p>d}O_\eps(\tnpd).
 \end{equation}
\end{proposition}
\begin{proof}
 {\it Step 1.} We first separate the ball $B_{r_1}$ from the annulus $A:=A_{r_1,R_n}$ (recall \eqref{defannulus}). We set
 \[
 n_A:=\hat{\mu}_n(A)\qquad \textrm{and then  } \kappa_A:=\frac{n_A}{\rho_n(A)}
 \]
and estimate thanks to triangle inequality and subadditivity,
\begin{multline}\label{separball}
 W_p^p\lt(\kappa_1 \one_{B_{r_1}}\rho_n+\sum_{j=2}^{j_n}\frac{1}{\stilde_j}\int_{\bStilde^{d-1}}\kappa_j^z \one_{\Delta_j^z} \rho_n d\stilde(z), n\rho_n\rt)\\
 \les W_p^p\lt(\kappa_1 \one_{B_{r_1}}\rho_n+\kappa_A \one_{A} \rho_n, n\rho_n\rt)\\
 +W_{A}^p\lt(\sum_{j=2}^{j_n}\frac{1}{\stilde_j}\int_{\bStilde^{d-1}}\kappa_j^z \one_{\Delta_j^z} \rho_n d\stilde(z), \kappa_A \rho_n\rt).
\end{multline}
To estimate the first term we use Proposition \ref{prop:radial1} with $\lambda(r)= r^{d-1}\rho_n(r)$, which is of type $\gamma=1-q$, $n_0=\hat{\mu}_n(B_{r_1})$ and $n_1=n_A$ so that
\[
 n\lambda(I_i)\ges n\gg1.
\]
By \eqref{estimradialabstract},
\begin{equation}\label{globalball}
 \EE\lt[W_p^p\lt(\kappa_1 \one_{B_{r_1}}\rho_n+\kappa_A \one_{A} \rho_n, n\rho_n\rt)\rt]\les_\eps n^{1-\frac{p}{2}}=o_\eps(\tnpd).
\end{equation}
In the next two steps we estimate the second right-hand side term in \eqref{separball}. We claim that 
\begin{equation}\label{claimgloballA}
 \EE\lt[W_{A}^p\lt(\sum_{j=2}^{j_n}\frac{1}{\stilde_j}\int_{\bStilde^{d-1}}\kappa_j^z \one_{\Delta_j^z} \rho_n d\stilde(z), \kappa_A \rho_n\rt)\rt]=o_\eps(\tnpd)+\one_{p>d}O_\eps(\tnpd).
\end{equation}
Combined with \eqref{separball} and \eqref{globalball}, this would conclude the proof of \eqref{eq:global}.\\

{\it Step 2.} The case $d\ge 3$ or $d=2$ and $p\ge 2$. Let us recall that by \eqref{defdelta}, $\delta_j=\delta$ for every $j\ge 2$ in this case. Set $\bar\sigma:=\tilde\sigma_j=\int_{\mathbb{\tilde S}^{d-1}}\mathbbm{1}_{D_{\delta}(z)}(e_1)d\stilde(z)$. Here, we first apply the angular estimate from Proposition \ref{prop:angle} to pass from $A$ to the angular sectors $\Lambda^z:=\Lambda_{r_1,R_n}^{z,\delta}$ and then decompose each angular sector using the radial estimate from Proposition \ref{prop:radial2}.
We thus define
\begin{equation}\label{defkappaz}
 n_z:=\hat{\mu}_n(\Lambda^z)\qquad \textrm{and then  } \qquad \kappa_z:=\frac{n_z}{\rho_n(\Lambda^z)}.
\end{equation}
Using triangle inequality we estimate
\begin{multline}\label{separglobalgeneric}
 W_{A}^p\lt(\sum_{j=2}^{j_n}\frac{1}{\sigma_j}\int_{\bStilde^{d-1}}\kappa_j^z \one_{\Delta_j^z} \rho_n d\stilde(z), \kappa_A  \rho_n\rt)
 \les W_A^p\lt(\frac{1}{\bar{\sigma}}\int_{\bStilde^{d-1}}\kappa_z \one_{\Lambda^z} \rho_n d\stilde(z), \kappa_A \rho_n\rt)\\
 + W_A^p\lt(\sum_{j=2}^{j_n}\frac{1}{\bar{\sigma}}\int_{\bStilde^{d-1}}\kappa_j^z \one_{\Delta_j^z} \rho_n d\stilde(z), \frac{1}{\bar{\sigma}}\int_{\bStilde^{d-1}}\kappa_z \one_{\Lambda^z} \rho_n d\stilde(z)\rt).
\end{multline}

{\it Step 2.1.} To estimate the first term in \eqref{separglobalgeneric}, we apply Proposition \ref{prop:angle} with $\lambda= (\rho_n/\rho_n(A))\one_A$ and $n=n_A$ (notice that $\EE[n_A]= n\rho_n(A)\simeq n$ and in this case $M_p\les 1$) to get
\begin{multline}\label{anggener}
 \EE\lt[W_A^p\lt(\frac{1}{\bar{\sigma}}\int_{\bStilde^{d-1}}\kappa_z \one_{\Lambda^z} \rho_n d\stilde(z), \kappa_A \rho_n\rt)\rt]\\
 \les \EE[n_A]^{1-\frac{p}{d-1}}\one_{d\geq 4}+\EE[n_A]^{1-\frac{p}{2}}[\one_{d=2}+(\log \EE[n_A])^{\frac{p}{2}}\one_{d=3}]\\
 \les \eta_n^{p,d-1}\stackrel{\eqref{deftaun}}{=} o(\tnpd).
\end{multline}

{\it Step 2.2.} To estimate the second term in \eqref{separglobalgeneric}, we first notice that for every $z\in \bStilde^{d-1}$,
\[
 \sum_j \kappa_j^z  \rho_n(\Dzj)\stackrel{\eqref{defkappajz}}{=} \hat{\mu}_n(\Lambda^z)\stackrel{\eqref{defkappaz}}{=} \kappa_z \rho_n(\Lambda^z)
\]
so that we may apply Lemma \ref{subad} and get
\begin{multline}\label{subadStep2.2}
 W_A^p\lt(\sum_{j=2}^{j_n}\frac{1}{\bar{\sigma}}\int_{\bStilde^{d-1}}\kappa_j^z \one_{\Delta_j^z} \rho_n d\stilde(z), \frac{1}{\bar{\sigma}}\int_{\bStilde^{d-1}}\kappa_z \one_{\Lambda^z} \rho_n d\stilde(z)\rt)\\
 \le \frac{1}{\bar{\sigma}}\int_{\bStilde^{d-1}} W_{\Lambda^z}^p\lt(\sum_{j=2}^{j_n}\kappa_j^z \one_{\Delta_j^z} \rho_n ,\kappa_z  \rho_n \rt)d\stilde(z).
\end{multline}
Fix $z\in \bStilde^{d-1}$. Let
\begin{equation}\label{defI}
 I:=(r_1,R_n), \qquad I_j:=(r_{j-1},r_j)
\end{equation}
and then
\begin{equation}\label{deflambda}
 \lambda(r):= \frac{\bar{\sigma}}{\rho_n(\Lambda^z)}r^{d-1} \rho_n(r) \one_{I},  \qquad  n_j:=\hat{\mu}_n(\Dzj), \qquad \textrm{and } \qquad \kappa_j:=\frac{n_j}{\lambda(I_j)}.
\end{equation}
By \eqref{defkappajz} and \eqref{defkappaz}, we have
\[
 \bar{\sigma} \kappa_j^z r^{d-1} \rho_n(r)=\kappa_j \lambda(r) \qquad \textrm{and } \qquad \bar{\sigma} \kappa_z  r^{d-1} \rho_n(r)=n_z \lambda(r)
\]
so that by radial symmetry of $\rho_n$,
\[
 W_{\Lambda^z}^p\lt(\sum_{j=2}^{j_n}\kappa_j^z \one_{\Delta_j^z} \rho_n ,\kappa_z  \rho_n \rt)=W^p_p\lt(\sum_{j=2}^{j_n} \kappa_j \one_{I_j} \lambda, n_z \lambda\rt).
\]
Let us notice that $n_z$ is a binomial random variable with
\[
 \EE\lt[n_z\rt]= n\rho_n(\Lambda^z)\simeq n \delta^{d-1}\gg 1
\]
so that with overwhelming probability, $n_z\ge  n \rho_n(\Lambda^z)/2$. In that event, for every $j\in[2,j_n]$,
\[
 n_z\lambda(I_j)\ges n \rho_n(\Dzj)\stackrel{\eqref{Enjz}}{\ges_\eps} 1.
\]
Therefore, since $\lambda$ is easily seen to be of  type $\gamma=1-q$ with constants depending only on $\eps$ and since when conditioned on $n_z$, $(n_j)_{j\in[2,j_n]}$ are  multinomial with parameters $(n_z\lambda(I_j))_{j\in[2,j_n]}$, we may apply Proposition \ref{prop:radial2} to get by definition \eqref{deflambda} of $\lambda$,
\begin{equation*}
\begin{split}
 \EE\lt[W_{\Lambda^z}^p\lt(\sum_{j=2}^{j_n}\kappa_j^z \one_{\Delta_j^z} \rho_n ,\kappa_z  \rho_n \rt)\rt] & \les_\eps \EE\lt[n_z\rt]^{1-\frac{p}{2}}\int_I r^{\frac{p(1-q)}{2}} \lambda(r)^{1-\frac{p}{2}} dr\\
&  \les_\eps (n\delta^{d-1})^{1-\frac{p}{2}}\begin{cases}
1 &\textrm{if } p<2\\
1+ R_n^{2-q}+\one_{q=2}\log R_n &\textrm{if } p=2\\
(R_n^{d+q-2} \rho(R_n))^{1-\frac{p}{2}} R_n^{2(1-q)} & \textrm{otherwise}.                                                                                                                                                                                                                                            \end{cases}
\end{split}
\end{equation*}
Plugging this back into \eqref{subadStep2.2} we find 
\begin{multline}\label{conclStep2.2}
 \EE\lt[W_A^p\lt(\sum_{j=2}^{j_n}\frac{1}{\bar{\sigma}}\int_{\bStilde^{d-1}}\kappa_j^z \one_{\Delta_j^z} \rho_n d\stilde(z), \frac{1}{\bar{\sigma}}\int_{\bStilde^{d-1}}\kappa_z \one_{\Lambda^z} \rho_n d\stilde(z)\rt)\rt]\\
 \les_\eps n^{1-\frac{p}{2}} \delta^{-(d-1)\frac{p}{2}} \begin{cases}
1 &\textrm{if } p<2\\
1+ R_n^{2-q}+\one_{q=2}\log R_n &\textrm{if } p=2\\
(R_n^{d+q-2} \rho(R_n))^{1-\frac{p}{2}}R_n^{2(1-q)} & \textrm{otherwise}.
\end{cases}
\end{multline}
Let us check that for $d\ge 3$ or $p\ge d$, the right-hand side is $o_\eps(\tnpd)+\one_{p>d}O_\eps(\tnpd)$. Combining this with \eqref{separglobalgeneric} and \eqref{anggener} would prove \eqref{claimgloballA} in this case. If $1\le p\le 2<d$, since $\delta$ is logarithmically small in $n$ and $n^{1-\frac{p}{d}}=\tnpd\gg n^{1-\frac{p}{2}}$, the statement is clear. If $p=d=2$, the right-hand side of \eqref{conclStep2.2} is of the order of
\[
 (\log n) (1+ (\log n)^{\frac{2}{q}-1}+ \log \log n)\ll (\log n)^{1+\frac{2}{q}}=\tnpd.
\]
We may thus focus on the case $p>2$. Recalling that $\delta\simeq \eps R_n^{-q}$, \eqref{condRn} and the definition \eqref{defqalpha} of $ \alpha, \beta$, there exists $\gamma\in \R$ such that  the right-hand side of \eqref{conclStep2.2} is of the order of
\[
 R_n^{(d-1)q-p(q-1)}(R_n^{-d(q-1)} n\rho(R_n))^{1-\frac{p}{2}}\les_\eps\begin{cases}
                                                     (\log n)^{\gamma}n^{-(1-\beta)(\frac{p}{2}-1)}& \textrm{if } 2<p<d\\
                                                     (\log n)^{\frac{d}{q}-1-\frac{1}{q}(\alpha-d(q-1))(\frac{d}{2}-1)}& \textrm{if } p=d\geq3\\(\log n)^{d-1-p(1-\frac{1}{q})}&\textrm{if } p>d.
                                                    \end{cases}
\]
Using that when $2<p<d$, $\beta<1$ and for $p=d$, $\alpha>d(q-1)$, and recalling the definition \eqref{deftaun} of $\tnpd$, the claim follows.\\

{\it Step 3.} We finally consider the case $d=2$ and $p<2$. In this case we will first split the annulus $A$ into the annuli $A_j:=A_{r_{j-1},r_j}$ and then split each annulus into the corresponding $\Dzj$. We introduce the notation
\begin{equation}\label{njStep3}
 n_j:=\hat{\mu}_n(A_j)\qquad \textrm{and then  } \qquad \kappa_{A_j}:=\frac{n_j}{\rho_n(A_j)}.
\end{equation}
By triangle inequality and subadditivity,
\begin{multline}\label{separglobalspecific}
 W_{A}^p\lt(\sum_{j=2}^{j_n}\frac{1}{\stilde_j}\int_{\bStilde^{d-1}}\kappa_j^z \one_{\Delta_j^z} \rho_n d\stilde(z), \kappa_A  \rho_n\rt)
 \\
 \les 
 W_A^p\lt(\sum_{j=2}^{j_n} \kappa_{A_j} \one_{A_j}\rho_n, \kappa_A \rho_n\rt)+\sum_{j=2}^{j_n} W_{A_j}^p\lt( \frac{1}{\stilde_j}\int_{\bStilde^{d-1}}\kappa_j^z \one_{\Delta_j^z} \rho_n d\stilde(z), \kappa_{A_j} \rho_n\rt).
\end{multline}

{\it Step 3.1.} In this step we estimate the first right-hand side term of \eqref{separglobalspecific}. We recall \eqref{defI} and set
\[
 \lambda(r):= \frac{1}{\rho_n(A)}r^{d-1} \rho_n(r) \one_{I} \qquad \textrm{and } \qquad \kappa_j:=\frac{n_j}{\lambda(I_j)}=\kappa_{A_j}\rho_n(A).
\]
We then have by radial symmetry of $\rho_n$,
\[
 W_A^p\lt(\sum_{j=2}^{j_n} \kappa_{A_j} \one_{A_j}\rho_n, \kappa_A \rho_n\rt)=W_p^p\lt(\sum_{j=2}^{j_n} \kappa_j \one_{I_j} \lambda, n_A \lambda\rt).
\]
As above we notice that with overwhelming probability, 
\[
 n_A\ge \EE[n_A]/2\ges n
\]
and in this event, for every $j\ge 2$ (recall that $d=2$ and that $r\mapsto r^{2-q} \rho(r)$ is decreasing in $[r_1,+\infty)$),
\begin{multline*}
 n_A\lambda(I_j)\ges n \rho_n(A_j)\ges_\eps n r_j (r_j-r_{j-1})\rho(r_j)\stackrel{\eqref{defrj}\&\eqref{defdelta}}{\ges_\eps} n r_j^{2-q} \rho(r_j)\\
 \ge n R_n^{2-q} \rho(R_n)\stackrel{\eqref{condRn}\&\eqref{defbarRn}}{\ges} (\log n)^{\frac{2}{q}-1}n^{1-\beta}\gg1.
\end{multline*}
Since $\lambda$ is of type $\gamma=1-q$ with constants depending only on $\eps$, we find from Proposition \ref{prop:radial2},
\begin{equation}\label{radialspecific}
 \EE\lt[W_A^p\lt(\sum_{j=2}^{j_n} \kappa_{A_j} \one_{A_j}\rho_n, \kappa_A \rho_n\rt)\rt]\les_\eps \EE[n_A]^{1-\frac{p}{2}} \int_{I} r^{p\frac{(1-q)}{2}}\lambda(r)^{1-\frac{p}{2}}dr
 \les n^{1-\frac{p}{2}}=o_\eps(\tnpd).
\end{equation}

{\it Step 3.2.} In this final step we estimate the second right-hand side term in \eqref{separglobalspecific}. For fixed $j$ we apply Proposition \ref{prop:angle} with $\delta=\pi$, $\lambda=\rho_n/\rho_n(A_j) \one_{A_j}$ and $n=n_j$ to get (recall that $d=2$)
\begin{multline*}
 \EE\lt[W_{A_j}^p\lt( \frac{1}{\stilde_j}\int_{\bStilde^{d-1}}\kappa_j^z \one_{\Delta_j^z} \rho_n d\stilde(z), \kappa_{A_j} \rho_n\rt)\rt]\les \EE[n_j]^{1-\frac{p}{2}}\int_{A_j} |x|^p d\lambda
 \les r_j^p (n\rho_n(A_j))^{1-\frac{p}{2}}.
\end{multline*}
Writing that by \eqref{defrj} and \eqref{defdelta},
\[
 r_j^p (n\rho_n(A_j))^{1-\frac{p}{2}}\simeq_\eps r_j^{1+\frac{pq}{2}} (n\rho(r_j))^{1-\frac{p}{2}} (r_j-r_{j-1})
\]
we have after summation,
\begin{multline}\label{Step3.2}
 \EE\lt[\sum_{j=2}^{j_n} W_{A_j}^p\lt( \frac{1}{\stilde_j}\int_{\bStilde^{d-1}}\kappa_j^z \one_{\Delta_j^z} \rho_n d\stilde(z), \kappa_{A_j} \rho_n\rt)\rt]\les_\eps n^{1-\frac{p}{2}} \sum_{j=2}^{j_n} r_j^{1+\frac{pq}{2}} \rho(r_j)^{1-\frac{p}{2}} (r_j-r_{j-1})\\
 \les_\eps n^{1-\frac{p}{2}} \int_{r_1}^{R_n} r^{1+\frac{pq}{2}} \rho(r)^{1-\frac{p}{2}}dr \les n^{1-\frac{p}{2}}=o_\eps(\tnpd).
\end{multline}
Combining \eqref{separglobalspecific}, \eqref{radialspecific} and \eqref{Step3.2} together concludes the proof of \eqref{claimgloballA} also in this case.

\end{proof}
\begin{remark}\label{rem:diff}
 Let us point out that if we tried to use the argument from {\it Step 2} in the case $1\le p<2$ and $d=2$ we would have an issue with the error coming from the radial part of the transport. Indeed, the error term coming from \eqref{conclStep2.2} would be much larger than $\tnpd=n^{1-p/2}(\log n)^{p/2}$. Moreover, if we kept $\delta_j=\delta$ we would have also an issue with the term \eqref{Step3.2}. Vice-versa, if we used the argument from {\it Step 3} in the case $p>d$ we would have an issue with the error coming from the angular part of the transport when $d=2$ or $d=3$ (see {\it Step 3.2} and \eqref{eq:angle}). In the remaining cases, both arguments give similar conclusions. 
\end{remark}

\section{Concentration inequalities}
In this section we briefly prove some concentration inequalities.
For $p\le d$ and $\eps>0$, let
\begin{equation}\label{defxin}
 \xi_n:=\begin{cases}
        n^{\frac{1}{p}-\frac{1}{2}} \tnpd^{-\frac{1}{p}} & \textrm{if } p\le 2\\
        \tnpd^{\frac{\eps-1}{2}}\one_{p<d} + \tnpd^{-\frac{1}{d}+\eps}\one_{p=d} &\textrm{otherwise}.
       \end{cases}
\end{equation}
Notice that $\xi_n\to 0$ as $n\to \infty$.
\begin{theorem}\label{thm:concentration}
Let $\{X_i\}_{i=1}^{\infty}$ be iid random variables with common distribution $\rho$ satisfying a Poincar\'e inequality and $\mu_n:=\sum_{i=1}^n\delta_{X_i}$. Then for every $1\leq p\leq d$, $r\ge 2$ and $t>0$, we have
\begin{eqnarray*}
\PP\lt( \tnpd^{-1}|W_p^p(\mu_n,n\rho)-\EE\lt[W_p^p(\mu_n,n\rho)\rt]|>t\rt)^{\frac{1}{r}}\les_{r,\eps} t^{-1} \xi_n.
\end{eqnarray*}
In particular we have complete (hence $\PP$ a.s.) convergence for $p<d$.
\end{theorem}
\begin{proof}
By Chebyshev inequality it is enough to prove that for every $r\ge 2$,
\begin{equation}\label{claimmomentbound}
 \mathbb{E}[|W_p^p(\mu_n,n\rho)-\mathbb{E}[W_p^p(\mu_n,n\rho)]|^r]^{\frac{1}{r}}\les_{\eps,r} \tnpd \xi_n.
\end{equation}

We start by recording the following a result from \cite[Proposition 5.3]{goldman2021convergence}: if $\pi$ is the optimal transport plan between $\mu_n$ and $n\rho$  and $r\ge 2$, we have
\begin{eqnarray}\label{conc1}
\mathbb{E}[|W_p^p(\mu_n,n\rho)-\mathbb{E}[W_p^p(\mu_n,n\rho)]|^r]^{\frac{1}{r}}\lesssim\mathbb{E}\left[\lt(\int_{\mathbb{R}^d\times\mathbb{R}^d}|x-y|^{2(p-1)}d\pi(x,y)\rt)^{\frac{r}{2}}\right]^{\frac{1}{r}}.
\end{eqnarray}
{\it Step 1.} The case $p\le 2$. In this case, since $2(p-1)\le p$, we may use H\"older inequality to upgrade \eqref{conc1} into
\begin{eqnarray}\label{conc1bis}
\mathbb{E}[|W_p^p(\mu_n,n\rho)-\mathbb{E}[W_p^p(\mu_n,n\rho)]|^r]^{\frac{1}{r}}\lesssim n^{\frac{1}{p}-\frac{1}{2}}\mathbb{E}\left[\lt(W_p^p(\mu_n,n\rho)\rt)^{\frac{r(p-1)}{p}}\right]^{\frac{1}{r}}.
\end{eqnarray}
We now prove by induction that for every $r\in \N$ with $r\ge 2$, we have
\begin{eqnarray}\label{claimple2}
\mathbb{E}[|W_p^p(\mu_n,n\rho)-\mathbb{E}[W_p^p(\mu_n,n\rho)]|^r]^{\frac{1}{r}}\lesssim_r n^{\frac{1}{p}-\frac{1}{2}}\tnpd^{1-\frac{1}{p}}.
\end{eqnarray}
By definition \eqref{defxin} of $\xi_n$ this would conclude the proof of \eqref{claimmomentbound} in this case.
Before embarking on the proof, let us point out that as a consequence of \eqref{claimple2}, $\xi_n\to 0$ and triangle inequality we have
\begin{equation}\label{coroclaimple2}
 \mathbb{E}[\lt(W_p^p(\mu_n,n\rho)\rt)^r]^{\frac{1}{r}}\lesssim_r \tnpd.
\end{equation}
We start by proving \eqref{claimple2} in the case $r=2$. In this case, since $2(p-1)/p\le 1$ it is a direct consequence of \eqref{conc1bis}, H\"older inequality and $\EE[W_p^p(\mu_n,n\rho)]=\tnpd$ by definition.\\
We now assume that \eqref{claimple2} holds up to $r-1$. Let us show that it holds also for $r$. Since
\[
 r\frac{p-1}{p}\le \frac{r}{2}<r-1,
\]
We have by \eqref{coroclaimple2} with $r-1$ instead of $r$ and H\"older inequality,
\[
 \mathbb{E}\left[\lt(W_p^p(\mu_n,n\rho)\rt)^{\frac{r(p-1)}{p}}\right]^{\frac{1}{r}}\le \mathbb{E}\left[\lt(W_p^p(\mu_n,n\rho)\rt)^{r-1}\right]^{\frac{p-1}{p(r-1)}}\les_r \tnpd^{1-\frac{1}{p}}
\]
so that by \eqref{conc1bis}, we get that \eqref{claimple2} also holds for $r$.\\

\medskip

{\it Step 2.} The case $p>2$. Notice that in particular we are then  only concerned with $d\ge 3$. For  $0<\eta<1$ we write $2(p-1)=(1-\eta) p+\eta m$ where
\begin{equation}\label{defm}
 m=\frac{2(p-1)-(1-\eta) p}{\eta}.
\end{equation}
Using  triangle inequality in the form
\[
 |a-b|^m\le (|a|+|b|)^m\le 2^m \max(|a|^m,|b|^m)\le 2^m(|a|^m+|b|^m)
\]
we have
\[
 \EE\lt[\int_{\R^d\times\R^d}|x-y|^m d\pi\rt]\les  2^m n \int_{\R^d} |y|^m d\rho\les (2^qm)^{\frac{m}{q}} n.
\]
In the last line we used that if $\rho=\exp(-V)$ with $V$ satisfying \eqref{hypV}, then $(\int|x|^Md\rho)^{\frac{1}{M}}\simeq M^{1/q}$. Indeed, this follows from $\rho\les \exp(-C |x|^q)$ and a change of variables to reduce to the case $q=1$ where the computation is explicit.
Combining this with H\"older inequality we can derive from \eqref{conc1},
\begin{multline}\label{conc1ter}
\mathbb{E}[|W_p^p(\mu_n,n\rho)-\mathbb{E}[W_p^p(\mu_n,n\rho)]|^r]^{\frac{1}{r}}
\lesssim (2^qm)^{\frac{m\eta}{2 q}} n^{\frac{\eta}{2}} \mathbb{E}\left[\lt(W_p^p(\mu_n,n\rho)\rt)^{\frac{(1-\eta) r}{2}}\right]^{\frac{1}{r}}.
\end{multline}
We claim that for every $\eps\ll1$ and  $r\in \N$ with $r\ge 2$,
\begin{eqnarray}\label{claimpge2}
\mathbb{E}[|W_p^p(\mu_n,n\rho)-\mathbb{E}[W_p^p(\mu_n,n\rho)]|^r]^{\frac{1}{r}}\lesssim_{r,\eps} \tnpd^{\frac{1+\eps}{2}}\one_{p<d} +\tnpd^{1-\frac{1}{d}+\eps}\one_{p=d}.
\end{eqnarray}
As in the case $p\le 2$ this would conclude the proof of \eqref{claimmomentbound} by definition \eqref{defxin} of $\xi_n$. Notice that if \eqref{claimple2} holds for some $r$, then  \eqref{coroclaimple2}
also holds for the same $r$ by triangle inequality. We fix $1\gg\eps>0$, $0<c<2\eps d/q$ and choose
\[
 \eta:=\begin{cases}
        \eps\lt(\frac{d}{p}-1\rt) &\textrm{if } p<d\\
        c \frac{\log \log n}{\log n} &\textrm{if } p=d.
       \end{cases}
\]
Notice that this choice of $\eta$ is such that
\begin{equation}\label{choiceeta}
 (2^qm)^{\frac{m\eta}{2q}} n^{\frac{\eta}{2}} \tnpd^{\frac{1-\eta}{2}}\les_\eps \tnpd^{\frac{1+\eps}{2}}\one_{p<d} +\tnpd^{1-\frac{1}{d}+\eps}\one_{p=d}.
\end{equation}

We prove \eqref{claimpge2} by induction and start with $r=2$. In this case, from \eqref{conc1ter} and H\"older inequality we have
\begin{multline*}
 \mathbb{E}[|W_p^p(\mu_n,n\rho)-\mathbb{E}[W_p^p(\mu_n,n\rho)]|^r]^{\frac{1}{r}}
\lesssim (4m)^{\frac{m\eta}{4}} n^{\frac{\eta}{2}} \mathbb{E}\left[W_p^p(\mu_n,n\rho)\right]^{\frac{1-\eta}{2}}\\
=(4m)^{\frac{m\eta}{4}} n^{\frac{\eta}{2}} \tnpd^{\frac{1-\eta}{2}}\stackrel{\eqref{choiceeta}}{\les_\eps}\tnpd^{\frac{1+\eps}{2}}\one_{p<d} +\tnpd^{1-\frac{1}{d}+\eps}\one_{p=d}.
\end{multline*}
This proves \eqref{claimpge2} when $r=2$.\\

Assume now that \eqref{claimpge2} holds for $r-1$. We then have
\[
 \frac{(1-\eta)r}{2}\le \frac{r}{2}\le r-1
\]
so that by \eqref{coroclaimple2} with $r-1$ and H\"older we have
\[
 \mathbb{E}\left[\lt(W_p^p(\mu_n,n\rho)\rt)^{\frac{(1-\eta) r}{2}}\right]^{\frac{1}{r}}\le \mathbb{E}\left[\lt(W_p^p(\mu_n,n\rho)\rt)^{r-1}\right]^{\frac{1-\eta}{2(r-1)}}\les_r \tnpd^{\frac{1-\eta}{2}}.
\]
Plugging this into \eqref{conc1ter} and recalling \eqref{choiceeta} concludes the proof of \eqref{claimpge2} also for $r$.

\end{proof}
\begin{remark}
 The proof  above fails when $p>d$. This may be easily seen from the fact that we have a gap in the rate, that is $\tnpd=(\log n)^{d-1-\frac{p}{2}}$ (that does not coincide with the correct rate when $p=d$), and this prevents us from finding a suitable $\eta=\eta(n)$ as we did before.
\end{remark}
\begin{remark}\label{rem:Poinc}
On the one hand, all log-concave distributions satisfy a Poincar\'e inequality by \cite{bobkov1999isoperimetric}. On the other hand, if $\rho$ satisfies a Poincar\'e inequality then it must have exponentially decaying tails by \cite{bobkov1997poincare}. Therefore, Theorem \ref{thm:concentration} mostly applies for potentials satisfying \eqref{hypV} with $q\ge 1$.
\end{remark}

\appendix
\section{Bounds for the matching problem in bounded domains}\label{sec:appendix}
We start by recording that arguing verbatim as in \cite[Proposition B.1]{TreMar}, we have the following de-Poissonization result
\begin{proposition}\label{prop:depoi}
 Let $\Omega$ be a bounded and connected set, $\lambda$ be a probability measure on $\Omega$, $\{X_i\}_{i=1}^\infty$ be i.i.d. random variables of law $\lambda$ and for every $n$, let $N_n$ be Poisson random variables of parameter $n$. Letting $\mu_n=\sum_{i=1}^n \delta_{X_i}$ and then for $W_*^p\in \{W_{\Omega}^p,\Wb_{\Omega}^p\}$,
 \[
  f(n):=\EE[W_*^p(\mu_n,n\lambda)],
 \]
we have for every $\alpha,\beta\in \R$,
\[
 \liminf_{n\to \infty} n^{\alpha} (\log n)^\beta f(n)\ge \liminf_{n\to \infty} n^{\alpha} (\log n)^\beta \EE[f(N_n)]
\]
and 
\[
 \limsup_{n\to \infty} n^{\alpha} (\log n)^\beta f(n)\le \limsup_{n\to \infty} n^{\alpha} (\log n)^\beta \EE[f(N_n)].
\]

\end{proposition}

We start by discussing upper bounds. Let us recall that if $X_i$ are i.i.d. random variables uniformly distributed in $Q_1=(0,1)^d$ and if $\mu_n=\sum_{i=1}^n \delta_{X_i}$ we have set
\[
 \betsup=\limsup_{n\to \infty} \frac{1}{\etnpd}\EE\lt[W_{Q_1}^p(\mu_n,n)\rt]
\]
where
\[
 \etnpd=n^{1-\frac{p}{d}}(1+\one_{d=2} (\log n)^{\frac{p}{2}}).
\]
It is now well-established, see e.g. \cite[Section 4]{Le17} that $\betsup\in(0,\infty)$ for every $d\ge 2$ and every $p\ge 1$. Notice that actually by \cite{AmStTr16,BaBo,DeScSc13,goldman2021convergence}, the limsup is a limit when $p=d=2$ (with $\betsup=1/(4\pi)$) or when $d\ge 3$. Combining the proof of \cite{TreMar} where the case of uniform densities is treated together with \cite[Proposition 2.4]{ambrosio2022quadratic} we have (see also  \cite{BaBo,DeScSc13,GolTrecombi}).
\begin{theorem}\label{th:upperboundappendix}
 Let $\Omega\subset \R^d$ be a bounded and connected open set either $C^2$ or convex and let $\lambda$ be a H\"older continuous probability density on $\overline{\Omega}$ bounded away from $0$. Let $X_i$ be i.i.d. random variables  with common law $\lambda$. then
 \begin{equation}\label{upperboundappendix}
  \limsup_{n\to \infty}  \frac{1}{\etnpd}\EE\lt[W_{\Omega}^p(\mu_n,n\lambda)\rt]\le \betsup \int_{\Omega} \lambda^{1-\frac{p}{d}}.
 \end{equation}

\end{theorem}

We now turn to the lower bounds and  recall that
\[
 \betinf=\liminf_{n\to \infty} \frac{1}{\etnpd}\EE\lt[\Wb_{Q_1}^p(\mu_n,n)\rt].
\]
We first prove that $\betinf\in(0,\infty)$.
\begin{proposition}\label{prop:betainf}
 For every $d\ge 2$ and every $p\ge1$, we have $\betinf\in(0,\infty)$. Moreover, the liminf in the definition of $\betinf$ is a limit at least if $p=d=2$ or if $d\ge3$.
\end{proposition}
\begin{proof}
Notice first that  $\betinf\le \betsup<\infty$. The fact that the limit exists and is equal to $1/(4\pi)$ in the case $p=d=2$ was proven in \cite[Proposition 3.1]{ambrosio2022quadratic} so we focus on the remaining cases. By scaling, Proposition \ref{prop:depoi} and H\"older inequality, it is enough to prove the following Poisson version of the statement: let $\mu$ be a Poisson point process of intensity one on $\R^d$ and set 
\[
 f_{p,d}(L):=\EE\lt[\frac{1}{|Q_L|}\Wb_{Q_L}^p\lt(\mu, \kappa\rt)\rt]
\]
where $\kappa:=\mu(Q_L)/|Q_L|$
then 
\begin{equation}\label{lowerboundcube}
 f_{1,d}(L)\ges (\log L)^{\frac{1}{2}}
\end{equation}
and for $d\ge 3$, $p\ge 1$, (the existence of the limit is part of statement)
\begin{equation}\label{lowerboundcubedgreater3}
\lim_{L\to \infty} f_{p,d}(L)=\betinf>0.
\end{equation}
For \eqref{lowerboundcube} we use \cite[Lemma 2.7]{huesmann2021there} to infer that 
\begin{multline*}
 f_{1,d}(L)=\EE\lt[\frac{1}{L^2} \sup\lt\{ \int_{Q_L} \zeta d(\mu-\kappa) \ : |\nabla \zeta|\le 1, \zeta=0 \textrm{ on } \partial Q_L\rt\}\rt]\\
 \ge \EE\lt[\frac{1}{L^2} \sup\lt\{ \int_{Q_L} \zeta d(\mu-\kappa) \ : |\nabla \zeta|\le 1, \zeta=0 \textrm{ on } \partial Q_L \ \& \ \int_{Q_L} \zeta=0\rt\}\rt]\\
 =\EE\lt[\frac{1}{L^2} \sup\lt\{ \int_{Q_L} \zeta d\mu \ : |\nabla \zeta|\le 1, \zeta=0 \textrm{ on } \partial Q_L \ \& \ \int_{Q_L} \zeta=0\rt\}\rt]\\
 \ges (\log L)^{\frac{1}{2}}.
\end{multline*}
We now turn to \eqref{lowerboundcubedgreater3}. Since for $\eps\in (0,1)$,
\[
 \EE\lt[\frac{1}{|Q_L|}\Wb_{Q_L}^p\lt(\mu, \kappa\rt)\rt]\ge (1-\eps)\EE\lt[\frac{1}{|Q_L|}\Wb_{Q_L}^p\lt(\mu, 1\rt)\rt] -\frac{C}{\eps^{p-1}} \EE\lt[\frac{1}{|Q_L|}\Wb_{Q_L}^p\lt(1, \kappa\rt)\rt]
\]
and by Lemma \ref{change},
\begin{equation}\label{estimconst}
 \EE\lt[\frac{1}{|Q_L|}\Wb_{Q_L}^p\lt(1, \kappa\rt)\rt]\les L^p\EE[|1-\kappa|^p]\les L^{p(1-\frac{d}{2})},
\end{equation}
in order to prove \eqref{lowerboundcubedgreater3} it is enough to prove it for
\[
 \what{f}_{p,d}(L):=\EE\lt[\frac{1}{|Q_L|}\Wb_{Q_L}^p\lt(\mu, 1\rt)\rt].
\]
Let $L>0$ and $m\in \N$. We divide $Q_{mL}$ into $m^d$ cubes $Q_i$ of sidelength $L$. By superadditivity, see Lemma \ref{superad}, and stationarity of $\mu$,
\[
 \what{f}_{p,d}(mL)\ge \sum_{i=1}^{m^d} \EE\lt[\frac{1}{|Q_{mL}|}\Wb_{Q_i}^p\lt(\mu, 1\rt)\rt]
 =\sum_{i=1}^{m^d} \frac{|Q_L|}{|Q_{mL}|} \EE\lt[\frac{1}{|Q_{L}|}\Wb_{Q_L}^p\lt(\mu, 1\rt)\rt]=\what{f}_{p,d}(L).
\]
Since $\what{f}_{p,d}$ is a continuous function the statement follows.

\end{proof}
\begin{remark}
 Notice that with the notation from the proof of Proposition \ref{prop:betainf}, we have using the fact that $\what{f}_{p,d}(mL)$ is increasing in $m$, triangle inequality and \eqref{estimconst},
 \[
  f(L)\le \betinf +\frac{C}{L^{\frac{d-2}{2}}}.
 \]
See also \cite[Theorem 4.1]{goldman2021convergence} for a similar estimate regarding $\betsup$.
\end{remark}

We may now prove the counterpart of Theorem \ref{th:upperboundappendix}.
\begin{theorem}\label{th:lowerboundappendix}
 Let $\Omega\subset \R^d$ be a bounded open set  let $\lambda$ a a H\"older continuous probability density on $\overline{\Omega}$ bounded away from $0$. Let $X_i$ be i.i.d. random variables  with common law $\lambda$. then
 \begin{equation}\label{lowerboundappendix}
  \liminf_{n\to \infty}  \frac{1}{\etnpd}\EE\lt[\Wb_{\Omega}^p(\mu_n,n\lambda)\rt]\ge \betinf \int_{\Omega} \lambda^{1-\frac{p}{d}}.
 \end{equation}

\end{theorem}
\begin{proof}
Let $1\gg \eps>0$ be fixed. And for $z\in \Z^d$ let $Q_z:= \eps(z+Q_1)$. Setting
\[Z:=\{z\in \Z^d \ : Q_z\subset \Omega\} \qquad \textrm{and } \qquad \Omega_\eps:=\cup_{z\in Z} Q_z,\]
we have by superadditivity,  see Lemma \ref{superad},
\[
 \EE\lt[\Wb_{\Omega}^p(\mu_n,n\lambda)\rt]\ge \sum_{z\in Z} \EE\lt[\Wb_{Q_z}^p(\mu_n,n\lambda)\rt].
\]
We claim that for every $z\in Z$,
\begin{equation}\label{claim:lowerappen}
 \liminf_{n\to \infty}  \frac{1}{\etnpd}\EE\lt[\Wb_{Q_z}^p(\mu_n,n\lambda)\rt]\ge (1-\omega(\eps))\betinf \int_{Q_z} \lambda^{1-\frac{p}{d}}.
\end{equation}
Summing over $z$ and letting then $\eps\to 0$, this would conclude the proof of \eqref{lowerboundappendix}. We now prove \eqref{claim:lowerappen} and set for simplicity $Q:=Q_z$. Letting $N:=\mu_n(Q)$ and for $(\hat{X}_i)_{i=1}^\infty$ i.i.d. random variables of law $\hat{\lambda}:=\lambda\restr Q/\lambda(Q)$, $\hat{\mu}_N:=\sum_{i=1}^N \delta_{\hat{X}_i}$, $\mu_n\restr Q$ and $\hat{\mu}_N$ have the same law thus
\[
 \EE\lt[\Wb_{Q}^p(\mu_n,n\lambda)\rt]=\EE\lt[\EE[\Wb_{Q}^p(\hat{\mu}_k,n\lambda)|N=k]\rt].
\]
By \cite[Proposition 2.4]{ambrosio2022quadratic}, there is a transport map $T: Q\mapsto Q$ with $T_\sharp \hat{\lambda}=1/|Q|$ and such that $\Lip T\le 1+C\omega(\eps)$. Letting $\mu'_n:=\sum_{i=1}^n \delta_{T(\hat{X}_i)}$ we then have for every $k\in \N$,
\begin{multline*}
 \EE[\Wb_{Q}^p(\hat{\mu}_k,n\lambda)]\ge (1-\omega(\eps))\EE\lt[\Wb_{Q}^p\lt(\mu'_k,\frac{n\lambda(Q)}{|Q|}\rt)\rt]\\
 \ge (1-\omega(\eps))\EE\lt[\Wb_{Q}^p\lt(\mu'_k,\frac{k}{|Q|}\rt)\rt]-\frac{C}{\eps^{p-1}} \EE\lt[\Wb_{Q}^p\lt(\frac{k}{|Q|},\frac{n\lambda(Q)}{|Q|}\rt)\rt]\\
 \stackrel{\eqref{eq:change}}{\ge} (1-\omega(\eps))\EE\lt[\Wb_{Q}^p\lt(\mu'_k,\frac{k}{|Q|}\rt)\rt]-C_\eps \frac{|n\lambda(Q)-k|^p}{(n\lambda(Q))^{p-1}}.
\end{multline*}
Taking the expectation we find
\begin{multline*}
 \EE\lt[\Wb_{Q}^p(\mu_n,n\lambda)\rt]=(1-\omega(\eps))\EE\lt[\Wb_{Q}^p\lt(\mu'_N,\frac{N}{|Q|}\rt)\rt]-C_\eps \EE\lt[\frac{|n\lambda(Q)-N|^p}{(n\lambda(Q))^{p-1}}\rt]\\
 \ge (1-\omega(\eps))\betinf|Q|^{\frac{p}{d}}\EE\lt[\eta_N^{p,d}\rt]-C_\eps n^{1-\frac{p}{2}}
\end{multline*}
where in the last line we used the definition of $\betinf$ and the fact that $N$ is a binomial random variable with $\EE[N]=n\lambda(Q)$. Since $n^{1-\frac{p}{2}}=o(\etnpd)$ and
\[
 \lim_{n\to \infty} \EE\lt[\frac{\eta_{N}^{p,d}}{\etnpd}\rt]=\lambda(Q)^{1-\frac{p}{d}},
\]
we find
\[
 \liminf_{n\to \infty}\frac{1}{\etnpd}\EE\lt[\Wb_{Q_z}^p(\mu_n,n\lambda)\rt]\ge (1-\omega(\eps))\betinf |Q|^{\frac{p}{d}}\lambda(Q)^{1-\frac{p}{d}}\ge (1-\omega(\eps))\betinf\int_{Q} \lambda^{1-\frac{p}{d}}.
\]
This concludes the proof of \eqref{claim:lowerappen}.
\end{proof}

\section*{Acknowledgements \& Funding}
F.P., M.G.\ and D.T.\ acknowledge the project G24-202 ``Variational methods for geometric and optimal matching problems'' funded by Università Italo Francese.

E.C.\ and D.T.\ acknowledge partial support by PNRR - MUR project PE00000013 - ``FAIR - Future Artificial Intelligence Research''.

D.T.\ acknowledges the MUR Excellence Department Project awarded to the Department of Mathematics, University of Pisa, CUP I57G22000700001,  the HPC Italian National Centre for HPC, Big Data and Quantum Computing - CUP I53C22000690001, the PRIN 2022 Italian grant 2022WHZ5XH - ``understanding the LEarning process of QUantum Neural networks (LeQun)'', CUP J53D23003890006, the INdAM-GNAMPA project 2024 ``Tecniche analitiche e probabilistiche in informazione quantistica''.

\bibliographystyle{acm}
\bibliography{OT.bib}

\end{document}